\documentclass[12pt]{article}
\usepackage{amsfonts}
\usepackage{amsmath, graphicx, epsfig,psfrag, verbatim, mathtools}
\usepackage{eucal,  amssymb,amsmath,amsthm,graphicx, amsfonts, latexsym,xcolor}
\usepackage{tikz}
\usetikzlibrary{shapes}
\usepackage[utf8]{inputenc}
\newcommand{\floor}[1]{\lfloor #1 \rfloor}
\def\E{{\rm E}}
\def\P{{\rm P}}

\def\re{{\rm e}}
\def\var{{\rm var}}
\def\cov{{\rm cov}}

\newcommand{\eqd}{\stackrel{D}{=}}

\def\an{a^{(n)}}
\def\An{A^{(n)}}
\def\Cn{C^{(n)}}
\def\Fn{F_n}

\def\Jn{J^{(n)}}
\def\Rn{R^{(n)}}
\def\Sn{S^{(n)}}
\def\Tn{T^{(n)}}
\def\Un{U^{(n)}}

\def\Xn{X^{(n)}}

\def\Zn{Z^{(n)}}

\def\Znb{\bar{Z}^{(n)}}

\def\bm{\boldsymbol{m}}
\def\bM{\boldsymbol{M}}
\def\bS{\boldsymbol{S}}
\def\bZ{\boldsymbol{Z}}

\def\bell{\boldsymbol{\ell}}

\def\mcCnch{\check{\mathcal{C}}^{(n)}}
\def\mcN{\mathcal{N}}
\def\mcB{\mathcal{B}}
\def\mcC{\mathcal{C}}
\def\mcG{\mathcal{G}}
\def\mcH{\mathcal{H}}
\def\mcS{\mathcal{S}}
\def\mcW{\mathcal{W}}

\newcommand{\bmeta}{\mbox{\boldmath{$\eta$}}}

\newcommand{\hfigwidth}{7.1cm}

\def\En{\mathcal{E}^{(n)}}
\def\Ent{\tilde{\mathcal{E}}^{(n)}}
\def\Ench{\check{\mathcal{E}}^{(n)}}
\def\Gn{\mathcal{G}^{(n)}}
\def\Gnt{\tilde{\mathcal{G}}^{(n)}}
\def\Snmc{\mathcal{S}^{(n)}}
\def\Cnmc{\mathcal{C}^{(n)}}

\def\Znch{\check{Z}^{(n)}}
\def\Cnch{\check{C}^{(n)}}
\def\Rnch{\check{R}^{(n)}}
\def\mcCnch{\check{\mathcal{C}}^{(n)}}

\def\mcN{\mathcal{N}}
\def\mcB{\mathcal{B}}
\def\mcBn{\mathcal{B}^{(n)}}
\def\mcBnch{\check{\mathcal{B}}^{(n)}}

\def\chin{\chi^{(n)}}
\def\chint{\tilde{\chi}^{(n)}}

\def\rbn{R^{(n)}_{\bullet}}
\def\abn{A^{(n)}_{\bullet}}
\def\abnt{\tilde{A}^{(n)}_{\bullet}}
\def\abnh{\hat{A}^{(n)}_{\bullet}}

\def\barzrn{\bar{R}^{(n)}_{\bullet}}
\def\barzan{\bar{A}^{(n)}_{\bullet}}

\def\bartn{\bar{T}^{(n)}}

\def\convp{\stackrel{{\rm p}}{\longrightarrow}}
\def\convD{\stackrel{{\rm D}}{\longrightarrow}}

\begin{document} \parskip=5pt plus1pt minus1pt \parindent=0pt
\title{An epidemic model on a network having two group structures with tunable overlap}
\author
{Frank Ball\footnote{School of Mathematical Sciences, University of Nottingham, frank.ball@nottingham.ac.uk}, Tom Britton\footnote{Department of Mathematics, Stockholm University, tom.britton@math.su.se} and Peter Neal\footnote{School of Mathematical Sciences, University of Nottingham, peter.neal@nottingham.ac.uk}}
\date{\today}
\maketitle

\begin{abstract}

A network epidemic model is studied. The underlying social network has two different types of group structures, households and workplaces, such that each individual belongs to exactly one household and one workplace.  The random network is constructed such that a parameter $\theta$ controls the degree of overlap between the two group structures: $\theta=0$ corresponding to all household members belonging to the same workplace and $\theta=1$ to all household members belonging to distinct workplaces. On the network a stochastic SIR epidemic is defined, having an arbitrary but specified infectious period distribution, with global (community), household and workplace infectious contacts.
 The stochastic epidemic model is analysed as the population size $n\to\infty$ with the (asymptotic) probability, and size, of a major outbreak obtained. These results are proved in greater generality than existing results in the literature by allowing for any fixed $0 \leq \theta \leq 1$, a non-constant infectious period distribution, the presence or absence of global infection and potentially (asymptotically) infinite local outbreaks.

\end{abstract}

{\it MSC2020 subject classifications:} Primary 92D30, 60K35; secondary 60J80, 05C80, 91D30. 

\section{Introduction}\label{sec-Intro}

Epidemic models (as well as epidemic outbreaks) are highly affected by heterogeneities in the community in which the epidemic spreads. There has been considerable research into incorporating population structure into epidemic models. 
One of the simplest and most commonly studied population structures is where the population is partitioned into small mixing groups, often referred to as households (Becker and Dietz~\cite{BD}). Infectious individuals make both local, within their household, and global, within the community at large, infectious contacts. The probability an infective infects a given housemate is typically an order of magnitude higher than the probability they infect a randomly chosen member of the population who resides in a different household. The (asymptotic) probability of a major outbreak and the size of a major outbreak as the size of the population (number of households) tends to infinity can be obtained, see Ball {\it et al.}~\cite{BMST}.

The households model provides a starting point for incorporating realistic population structure into epidemic models for human populations. Human populations contain more intricate social structure which presents challenges in modelling (Ball {\it et al.}~\cite{Ball_etal15}). Households have a geographical location and the spatial location or neighbourhood to which a household belongs can be important in defining the mixing between individuals in different households, and hence the spread of a disease.  Increased realism can be obtained by using a hierarchical, three-level-mixing model (Britton {\it et al.}~\cite{BKO11}), in which the population is partitioned into neighbourhoods, which are partitioned into households;
with different rates of transmission between an infective and those individuals they share a household with, a neighbourhood with or are a general member of the  population.  The focus in Britton {\it et al.}~\cite{BKO11} is on inference; Ouboter et al.~\cite{OMT16} analyse the asymptotic probability and size of a major outbreak for such a three-level-mixing model, with neighbourhoods corresponding to schools, when the infectious period is constant.

As well as close ties with people they live with, individuals have regular contact  with individuals through work or school, and hence, increased opportunities to spread a disease with such people. This has led to models with multiple partitions of the population, for example, Andersson~\cite{And99} and Ball and Neal~\cite{BN02}. In the case of two partitions of the population the model is often referred to as the households-workplaces model. A common assumption, to aid mathematical analysis, is that the two partitions are constructed to have no overlap (see, for example, Ball and Neal~\cite{BN02}, Pellis et al.~\cite{Pellis09,Pellis11}, Kubasch~\cite{Kubasch23} and Bansaye et al.~\cite{BDKV24}). That is, no two individuals belong to the same household and the same workplace. This contrasts with the hierarchical model where all individuals in the same household belong to the same neighbourhood.  Indeed, Ouboter et al.~\cite{OMT16} prove that if the sizes of households and workplaces are both constant, and the infectious period is also constant, then the asymptotic probability and  size of a major outbreak are both greater in the model with no overlap than in the model with complete overlap.

We take workplace as a generic term to cover a mixing group structure outside the household, and to include schools and kindergartens. Therefore, the workplace will typically be a larger unit than the household. There are clearly situations in which substantial overlap occurs between the household and workplace. For example, it is often the case that siblings go to the same kindergarden or school, thus creating overlap between the two group structures. A second example, which is the main inspiration to the current work, is in Patwardhan et al.~\cite{P23} where university students from a group of (major) classes are allocated into dorms. Here dorms play the role of household and classes correspond to workplaces. Patwardhan et al.~\cite{P23} numerically analyse a closely related model and one of the questions they address is whether students should be allocated dorms according to which class they follow (large overlap  - in the limit the hierarchical model) or the dorm allocation should be done more randomly (small overlap - in the limit the standard households-workplaces model) in order to minimize effects of epidemic outbreaks.  A key finding is that allocation according to classes can not only significantly reduce the severity of an epidemic but it also enhances the effectiveness of vaccinating individuals, chosen uniformly at random, in reducing disease spread.

In the current paper we consider the households-workplaces model with a tunable parameter, $\theta$. The model is constructed by first allocating households to workplaces, so that the population structure resembles the neighbourhood model with neighbourhoods labelled by workplace. Then each individual decides independently to be a {\it mover} with probability  $\theta$. The movers change their workplace membership but not the household to which they belong. The movers are then randomly allocated to workplace spaces vacated by movers. The non-movers, who we term {\it remainers}, remain in the workplace and household to which they were originally assigned. In this way the number of individuals in each workplace does not change. A full description of the model is provided in Section \ref{sec-model} with an example of a construction of a population structure given in Figure \ref{fig:complex} in Section~\ref{subsubsec-complex}.  The special cases $\theta=0$ and $\theta=1$ yield complete overlap and no overlap, respectively.

We analyse the spread of an SIR (susceptible $\to$ infective $\to$ recovered) epidemic, having an arbitrary but specified infectious period distribution, on a network of households and workplaces constructed as above. We assume that the sizes of households and workplaces are both constant and consider the asymptotic regime in which the total population size $n\to\infty$, keeping all other parameters fixed.  We consider two scenarios: (a) in which there is no global infection (i.e. with a randomly chosen member of the population), so infection is spread entirely through household and workplace contacts; and (b) where there is global infection, in addition to household and workplace contacts.  The papers cited above are concerned with scenario (b), though Andersson~\cite{And99} also considers scenario (a).  We consider epidemics initiated by a single infective, chosen uniformly from the population, and prove that under each scenario that, as $n \to \infty$, the probability of a major outbreak (i.e. one that infects at least $\log n$ individuals) converges to a constant $\rho$ and conditional upon the occurrence of a major outbreak, the fraction of the population that gets infected converges in probability to a strictly positive constant $z$. For each scenario, we provide recipes for calculating the constants $\rho$, $z$ and a threshold parameter that determines whether or not $\rho>0$.

Whilst the special cases $\theta=0$ and $\theta=1$ have been considered in the literature, there has been no general analysis of the case $0 < \theta <1$, where there is partial overlap between the two group structures. Moreover, even for the cases $\theta=0$ and $\theta=1$, the only rigorous proofs of a law of large numbers for the size of a major outbreak assume that the infectious period is constant, making the numbers of global, household and workplace infectious contacts made by an infective independent thus simplifying the analysis.  The fact that households and workplaces (may) overlap disallows analysing small independent components, so instead we study {\it local susceptibility sets} and {\it local infectious clumps} in the spirit of Ball and Neal~\cite{BN02}. In  Ball and Neal~\cite{BN02} it was noted that the local infectious clumps (ignoring global infectious contacts) for the households-workplaces model, with $\theta=1$, can be infinite, but attention was focussed on the simpler case where the local infectious clumps are almost surely finite. The infectious clumps have certain \emph{complexes} as there base unit, where a complex is defined as individuals related to a given workplace: the remainers in the original construction, the movers in the original construction (who hence now belong to a different workplace) and the movers joining the workplace (and hence had a different workplace at the start).  For both scenarios (a) and (b), the proof of the limiting major outbreak probability uses branching process approximations.  For scenario (a), the proof of the law of large numbers for the size of a major outbreak adapts ideas from the proof of a corresponding result for epidemics on random intersection graphs in Ball et al.~\cite{BST14}, while for scenario (b) it uses an extension of the embedding technique of Scalia-Tomba~\cite{ScaliaTomba85, ScaliaTomba90}.  In both cases, the above-mentioned complexes are fundamental to the extensions, which are far from straightforward.

Beside different notation, our model is similar to the model in Patwardhan et al.~\cite{P23}, although global infectious contacts in the latter are achieved via a stochastic block model. In Patwardhan et al.~\cite{P23} the list of individuals is gone through sequentially, each time the individual decides to be a mover with probability $r$ and if so switches workplace with a randomly selected other individual. As a consequence, a given individual remains in their original household if it is not a mover \emph{and} it is not selected to switch workplace with another mover, an event which has asymptotic probability $(1-r){\rm e}^{-r}$ as $n \to \infty$, whereas individuals remain in their original workplace with asymptotic probability $1-\theta$ in our model. (In Patwardhan et al.~\cite{P23}, a given individual may return to their workplace after making one or more switches but that happens with probability zero under the above asymptotic regime.)  Patwardhan et al.~\cite{P23} consider the effect of $r$ on the time course and final size of the epidemic with particular focus on the cases $r=0$ $(\theta =0)$ and $r=1$ $(\theta=1)$. Two mean-field approximations for the epidemic are provided in the supplementary material of Patwardhan et al.~\cite{P23}. However, both approximations perform rather poorly both for the disease dynamics and the final size, unless $r=1$ and the epidemic is highly infectious, see Patwardhan et al.~\cite{P23} supplementary material, Figures S8 and S9. Moreover, Patwardhan et al.~\cite{P23} note that a very important remark is that their findings indicated above concerning allocation of students to dorms are apparent only from their numerical simulations and do not emerge from numerical integration of their mean-field approximations.  By contrast using the framework of  Ball and Neal~\cite{BN02} and dealing with the additional mathematical technicalities of allowing (asympotically) infinite local infectious clumps and overlapping groups, we are able to obtain exact results for the (asympotic) final size and  outbreak probability, which demonstrate clearly (see Figure~\ref{fig:finalsize3} in Section~\ref{sec-illustrations}) that disease spread can be increasing in $\theta$, i.e.~decreasing in overlap.

The remainder of the paper is structured as follows. In the first two short sections (\ref{sec-model} and \ref{sec-main}) we define the model and state our main results. In Section \ref{sec-illustrations} we illustrate our results and investigate them numerically.
In Sections  \ref{sec-derivations} and  \ref{sec-pig0} we prove the main results for scenarios (b), global infection present, and (a), no global infection, respectively. We end with a short discussion in Section \ref{sec-disc}.

\section{The epidemic model with two group structures and tunable overlap}\label{sec-model}

\subsection{The network model}
\label{sec-networkmodel}

Consider a population of size $n$ (assumed to be large) having a network structure in which each individual belongs to exactly one household and exactly one workplace. To simplify things we assume all households have size $h\ge 2$ and all workplaces have size $w$, where we make the additional assumption that $w=dh$ for some integer $d\ge 1$ (so workplaces are at least as large as and a multiple the size of households).  For convenience we consider a population partitioned into households and also into workplaces but the two groups structures could be other possible partitions, such as households and schools, as in Ouboter et al.~\cite{OMT16}, or dorms and classes, as in  Patwardhan et al.~\cite{P23}.

The parameter $\theta$ determining the amount of overlap (or rather the lack of overlap) is defined as follows. First we let each workplace consist of $d$ households, so all individuals of the same household are also in the same workplace implying a complete overlap. Then each individual independently decides whether to be a {\it mover} with probability $\theta$ or a {\it remainer} with the remaining probability $1-\theta$. The workplace spots of movers are made vacant, and all movers are distributed uniformly at random among the vacant spots in the different workplaces, thus creating new workplace constellations but all workplaces still have size $w=dh$. As a consequence, the larger $\theta$, the fewer pairs of household members who are also work colleagues.

In summary, the parameters of the network (beside the overall population size $n$ which must be a multiple of $w$) are: $h$, the household size, $w=dh$, the workplace size, and $\theta$ the probability of being a mover (the tunable overlap can hence be defined as $1-\theta$).

\subsection{The epidemic model on the network}
\label{sec-epidemicmodel}

The epidemic is a stochastic continuous-time SIR model in which each infected individual is infectious for an independent random duration $I$: the length of the infectious period (possibly preceded by a latent period having an arbitrary distribution). During the infectious period an individual makes three types of infectious contacts at times of independent Poisson processes: at rate $\beta_H$ the individuals make infectious contact with other members of their household, each time choosing which individual independently and uniformly at random, implying the rate to a specific individual equals $\beta_H/(h-1)$. Additionally, the infectious individual makes contact with workplace individuals at rate $\beta_W$, each time with an independent, uniformly selected individual of the same workplace. Finally, the infectious individual makes global contacts at rate $\beta_G$, each time with an individual chosen independently and uniformly from the entire population. An infectious contact with a susceptible individual results in the latter getting infected (such individuals starting their infectious contact processes), whereas contacts with already infected people have no effect. After their infectious period is over, the individual recovers and becomes immune thus playing no further role in the epidemic.

The epidemic is initiated by one randomly selected individual being externally infected and the remaining individuals being susceptible. The epidemic runs its course until there are no infectious individuals present in the community when it stops. Each individual is then either recovered or susceptible, and the random distribution of how the recovered and susceptibles are distributed in the community defines the final outcome. The final \emph{number} infected is denoted $Z$ and referred to as the final size.

For ease of analysis, we assume that each global contact is with an individual chosen uniformly from the entire population, including also the infector itself.  Thus a global contact may be with an individual in the infector's own household or workplace.  This assumption represents no loss of generality and has no impact on our asymptotic results.

\subsection{Parametrization}

Without loss of generality we assume $\E[I]=1$ (meaning that time is measured in units of the mean infectious period). Our focus lies on $h$ and $w$ being fixed and letting $n\to\infty$. If $w$ also becomes large with $n$ we essentially have a households epidemic model, and if both $h$ and $w$ grow with $n$, the model is very close to being homogeneously mixing.

The epidemic parameters are the three infectious contact rates $\beta_H,\ \beta_W$ and $\beta_G$, and the infectious period distribution $I$ (assumed to have mean $\E[I]=1$). 
Later we  reparameterize the three infection rates to $\beta=\beta_H+\beta_W+\beta_G$, $\pi_G=\beta_G/(\beta_H+\beta_W+\beta_G)$ and $\pi_{H|G^c}=\beta_H/(\beta_H+\beta_W)$, thus being the overall contact rate, what fraction of the overall contacts rate are global contacts, and what fraction of the remaining (local) contacts are in households, respectively.  

\subsection{Notation}
\label{notation}
We use $\convD$ and $\convp$ to denote convergence in distribution and convergence in probability, respectively, and $\eqd$ to denote equal in distribution.  For a random variable, $X$ say, taking values in the positive integers, $f_X$ denotes its probability-generating function (PGF), i.e.~$f_X(s)=\E[s^X]$ ($0 \le s \le 1$).  For any real-valued random variable $Y$, $\mu_Y$ denotes its mean and $\phi_Y$ denotes its Laplace transform, so
$\phi_Y(\nu)=\E[\exp(-\nu Y)]$ $(\nu \ge 0)$.  For a positive integer $n$ and $p \in [0,1]$, ${\rm Bin}(n,p)$ denotes a binomial random variable with $n$ trials and success probability $p$.  For $\lambda>0$, ${\rm Exp}(\lambda)$ denotes an exponential random variable with rate $\lambda$ and hence mean $\lambda^{-1}$.  For $x \in \mathbb{R}$, $\floor{x}$ denotes the usual floor function, so $\floor{x}$ is the greatest integer $\le x$.  For a non-negative random variable $X$, we write $Y\sim {\rm Po}(X)$ when $Y$ has a Mixed-Poisson distribution with parameter distributed as $X$, i.e.~when
$\P(Y=k) = \E[X^k {\rm e}^{-X} / k!]$ for $k=0,1,2,\dots$.

\section{Main Results}\label{sec-main}

Let $Z^{(n)}$ denote the final size of the epidemic in a community of size $n$, i.e.\ the total number getting infected throughout the course of the epidemic. Our first result concerns the probability of a major outbreak. Here we define a major outbreak as $Z^{(n)} > \log n$, but any sequence $h_n$ such that $h_n\to\infty$ and $h_n/n \to 0$ would do equally well.
We also define reproduction numbers $R_*$ and $R_L$, which determine if the asymptotic probability of a major outbreak is strictly positive for the cases $\pi_G>0$ and $\pi_G=0$, respectively.

${\newtheorem{theorem}{Theorem}[section]}$
\begin{theorem}\label{prop-pmajor}
For the epidemic model with two group structures and tunable overlap defined in Section \ref{sec-model},
\begin{equation}
\P(Z^{(n)}>\log n) \to \rho \quad \text{as } n\to \infty.
\label{maj-prob}
\end{equation}

When $\pi_G>0$, the constant $\rho$ (the major outbreak probability) is  given by $\rho=1-\xi$, where $\xi$ is the smallest solution in $[0,1]$ of $\phi_A(\beta_G(1-s))=s$ and the Laplace transform $\phi_A(\nu)$ is given by~\eqref{equ:phiA} in Section~\ref{subsec-outbreakprob}. The major outbreak probability $\rho$ is strictly positive if and only if $R_*>1$, where $R_*$ is given by~\eqref{equ:RstarS} in Section~\ref{subsubsec-susset}.

When $\pi_G=0$, the major outbreak probability $\rho$ is given by~\eqref{eq:pig0:extinct} in Section~\ref{sec-pig0-genI} and $\rho$ is strictly positive if and only if $R_L>1$, where $R_L$ is given by~\eqref{equ:zetaS} in Section~\ref{subsubsec-susset}.
\end{theorem}

Our next result concerns the limiting final size in case there is a major outbreak. 

${\newtheorem{theorem-finsize}[theorem]{Theorem}}$
\begin{theorem-finsize}\label{prop-finsize}
Consider the epidemic model with two group structures and tunable overlap defined in Section \ref{sec-model}. Suppose that $R_L>1$ if $\pi_G=0$, and that $R_*>1$ and ${\var}(I)<\infty$ if $\pi_G>0$.  We then have
\begin{equation}
\frac{Z^{(n)}}{n}|Z^{(n)}>\log n \convp z \quad \text{as } n\to \infty,
\label{LLN}
\end{equation}
where the constant $z$ (the limiting final size) is strictly positive.

If $\pi_G>0$, the limiting final size $z$ is given by the unique strictly positive root of the equation~\eqref{equ:z} in Section~\ref{sec-genframework}, with $f_S(s)$ being given by~\eqref{equ:progenypgfS} in Section \ref{subsubsec-susset}.

If $\pi_G=0$, the limiting final size $z$ is given by~\eqref{eq:pig0:extinctS} in Section~\ref{sec-pig0-clumpsus}.
\end{theorem-finsize}
${\newtheorem{remark}[theorem]{Remark}}$
\begin{remark}
It is worth noting that if the infectious period is constant, i.e.~$I \equiv 1$, then the major outbreak probability $\rho$ in~\eqref{maj-prob} is identical to the final size quantity $z$ in~\eqref{LLN}, but this does not hold when the infectious period is random. 
\end{remark}

In Section \ref{sec-illustrations} we compute $z$ (and hence also $\rho$ if $I \equiv 1$) numerically for various parameter settings and study how it depends on the different model parameters.

\section{Numerical illustrations}\label{sec-illustrations}
In this section we present numerical results which show the usefulness of the asymptotic results in Theorems~\ref{prop-pmajor} and~\ref{prop-finsize} for finite population size $n$, and illustrate the dependence of the limiting final size $z$ on the probability an individual is a mover, $\theta$, the number of households that initially comprise a workplace, $d=w/h$, and the distribution of the infectious period, $I$.

\begin{figure}
\begin{center}
\begin{tabular}{ccc}
 \includegraphics[width=\hfigwidth]{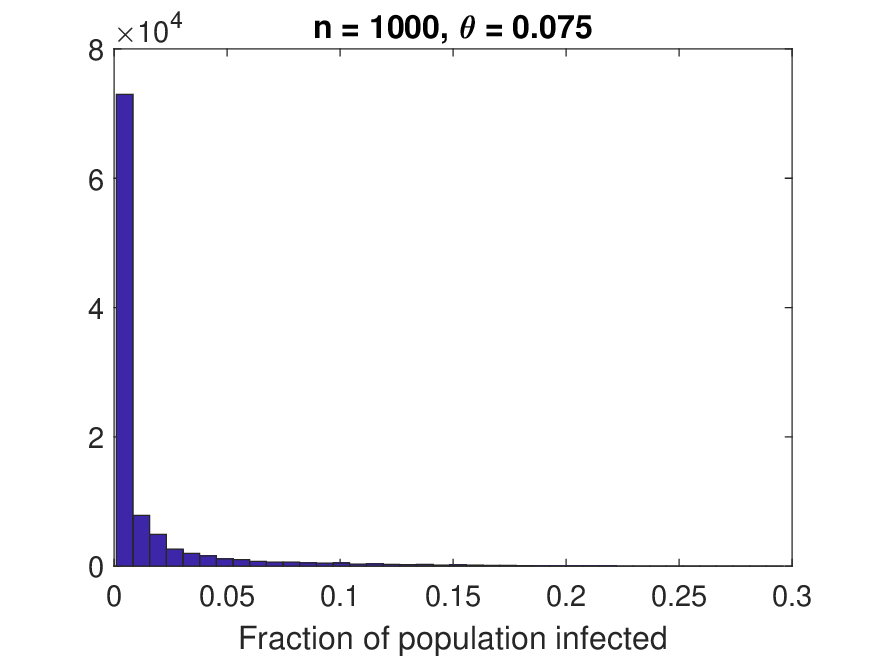} &
 \includegraphics[width=\hfigwidth]{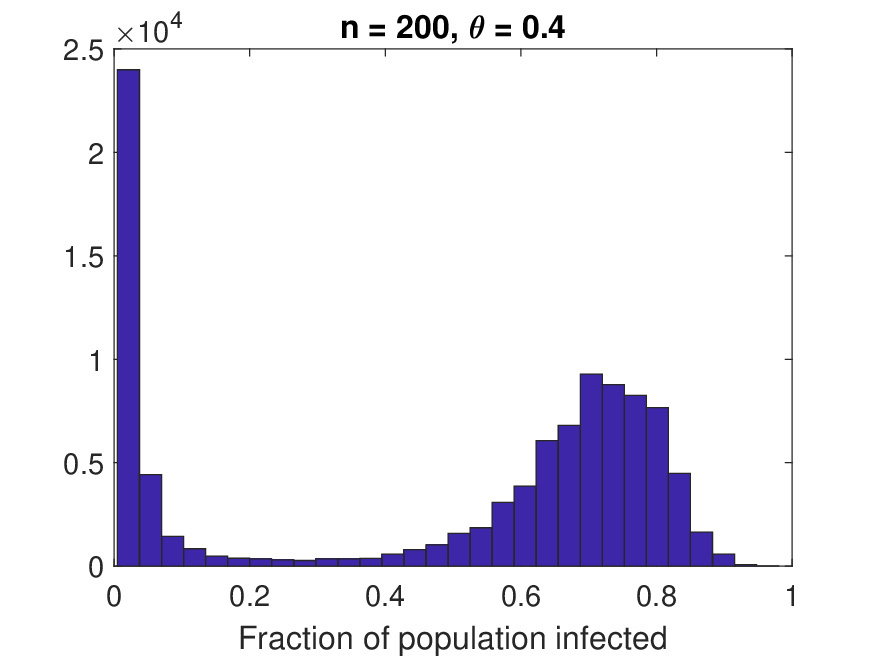} \vspace{.5 cm} \\
 \includegraphics[width=\hfigwidth]{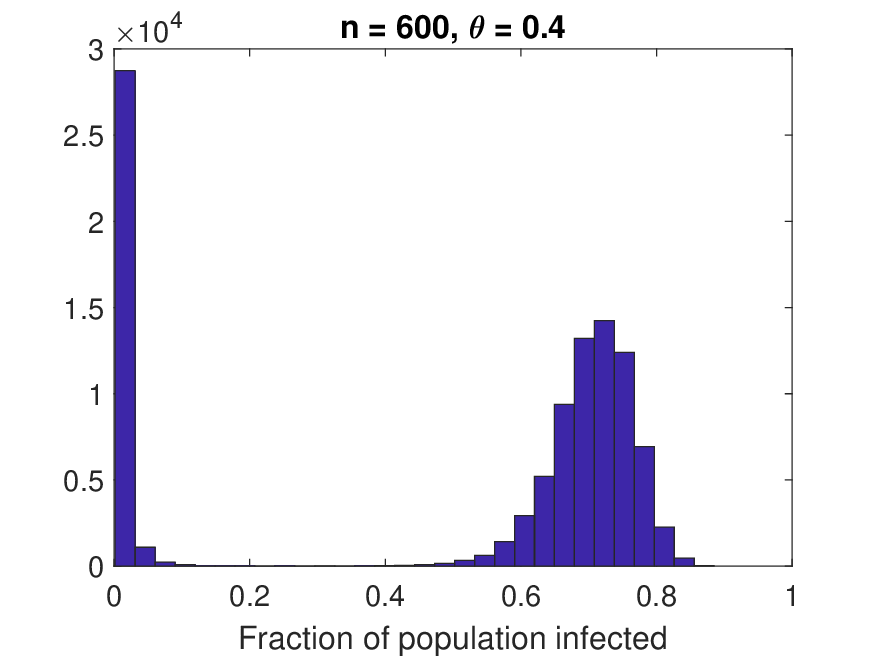} &
 \includegraphics[width=\hfigwidth]{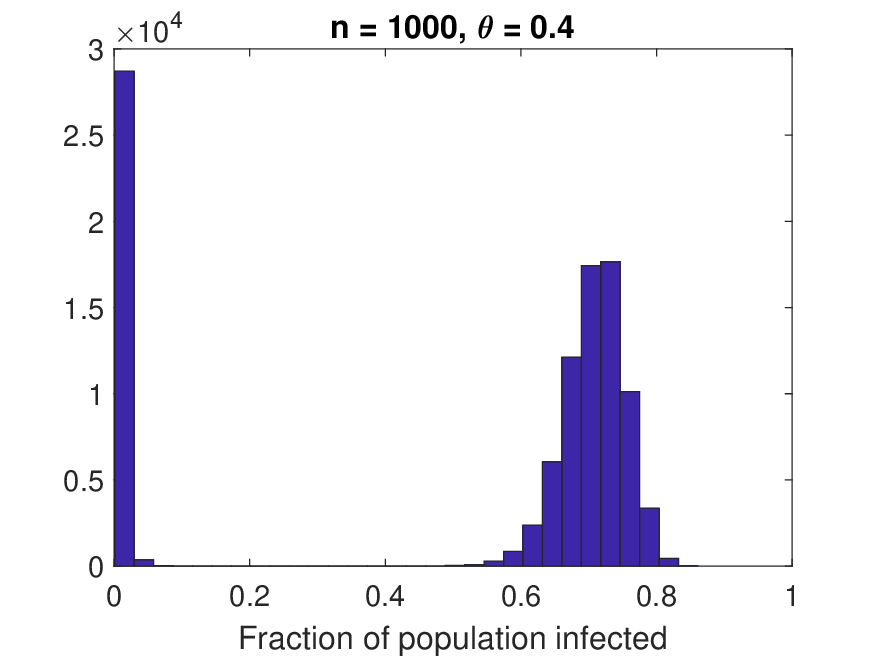}
\end{tabular}
\end{center}
\caption{Histograms of simulated fraction of the population infected for epidemics with $d=1$, $h=w=4$, $\beta=3$, $\pi_G=0.025$, $\pi_{H|G^c}=0.5$ and $I \equiv 1$. Each plot is based on 100,000 simulations. See text for further details.
}
\label{fig:finalsize1}
\end{figure}

In Figure~\ref{fig:finalsize1}, histograms are shown of the fraction of the population infected by epidemics with $d=1$, $h=w=4$, $\beta=3$, $\pi_G=0.025$, $\pi_{H|G^c}=0.5$ and $I \equiv 1$.  Each plot is based on 100,000 simulations of epidemics with one initial infective chosen uniformly at random from the population, with the remainder of the population being susceptible.  In the top-left plot, $\theta=0.075$ (large overlap) and $n=$ 1,000, so the population is partitioned into 250 households of size 4, and also into 250 workplaces of size 4.  With these parameter values, $R_*=0.6541$, so the epidemic is subcritical as is clear from the shape of the histogram.  In the other three plots, $\theta=0.4$ (smaller overlap), for which $R_*=\infty$ as the corresponding epidemic with only local infection (i.e.~$\beta_G=0$ and $\beta_H$ and $\beta_W$ are unchanged) is supercritical.   In each of these plots the final size distribution is bimodal, reflecting minor and major outbreaks.  When $n=$ 1,000, there is a sharp distinction between minor and major epidemics, with none of the 100,000 simulated epidemics having a final size in $(172, 358)$.  The distinction is less clear when $n=600$ and when $n=200$ the cut-off between minor and major epidemics is far from clear.

\begin{figure}
\begin{center}
\begin{tabular}{ccc}
$d=1$ & $d=1$ \\
 \includegraphics[width=\hfigwidth]{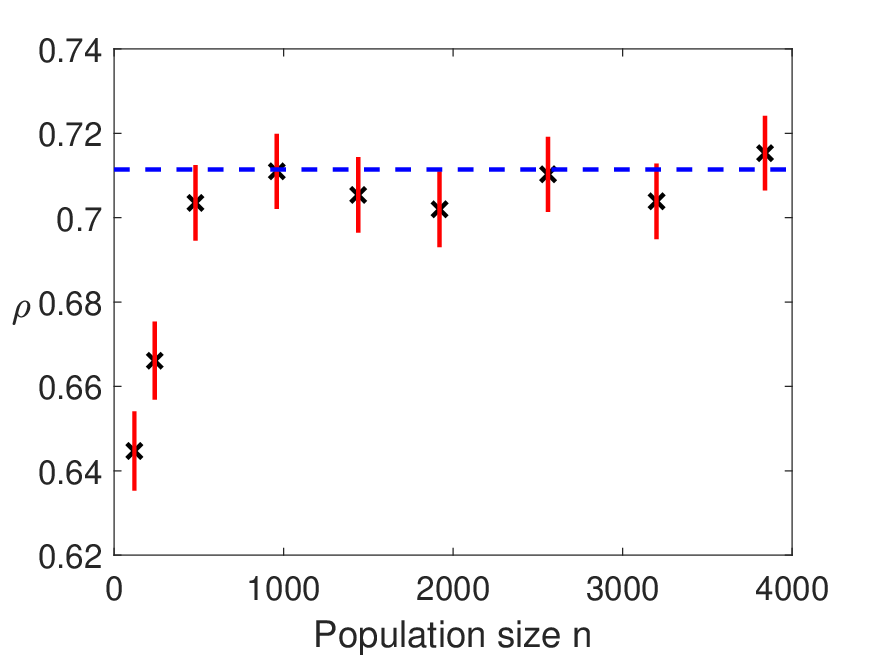} &
 \includegraphics[width=\hfigwidth]{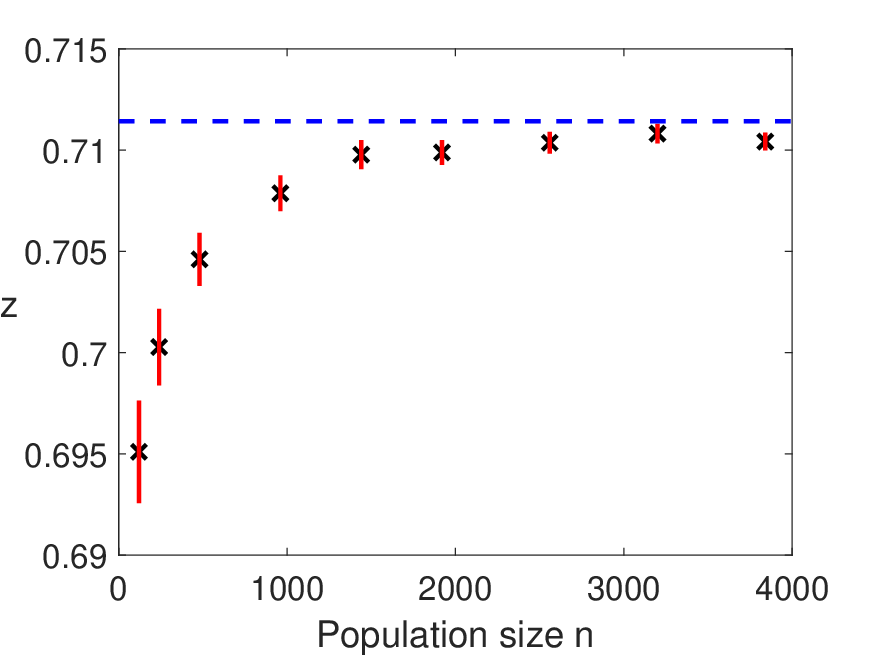} \vspace{.5 cm} \\

 $d=2$ & $d=2$ \\

 \includegraphics[width=\hfigwidth]{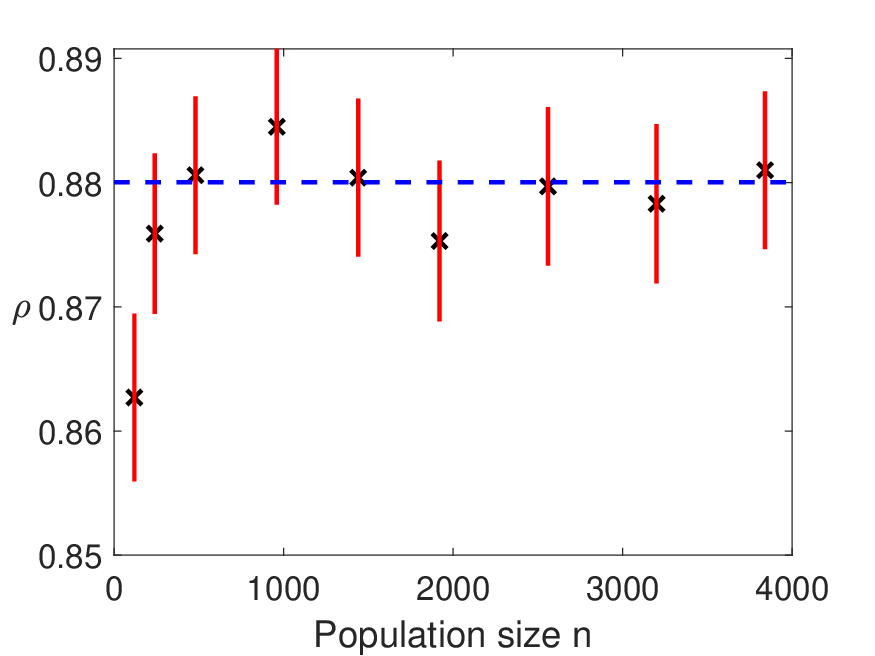} &
 \includegraphics[width=\hfigwidth]{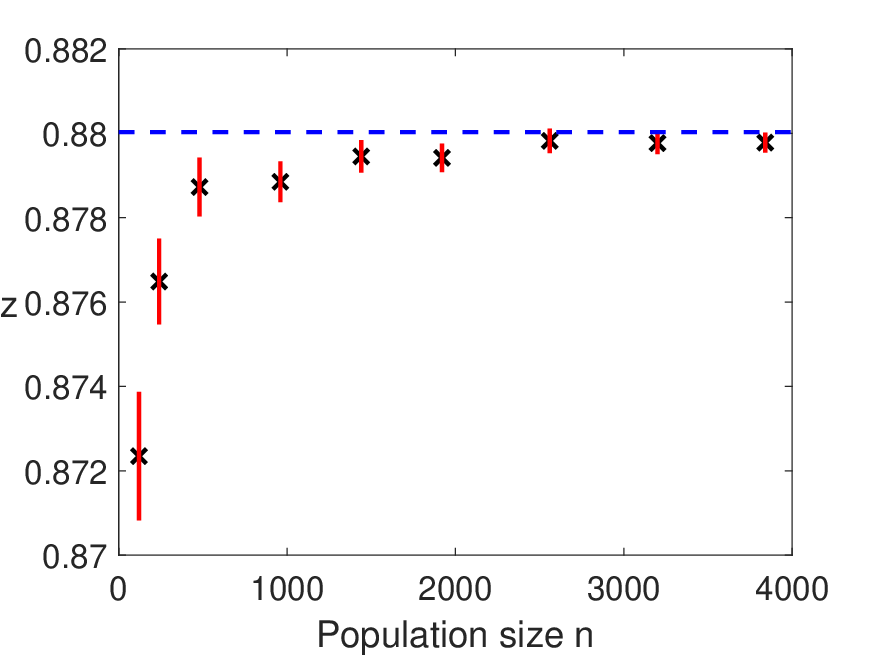} \vspace{.5 cm} \\

 $d=3$ & $d=3$ \\

 \includegraphics[width=\hfigwidth]{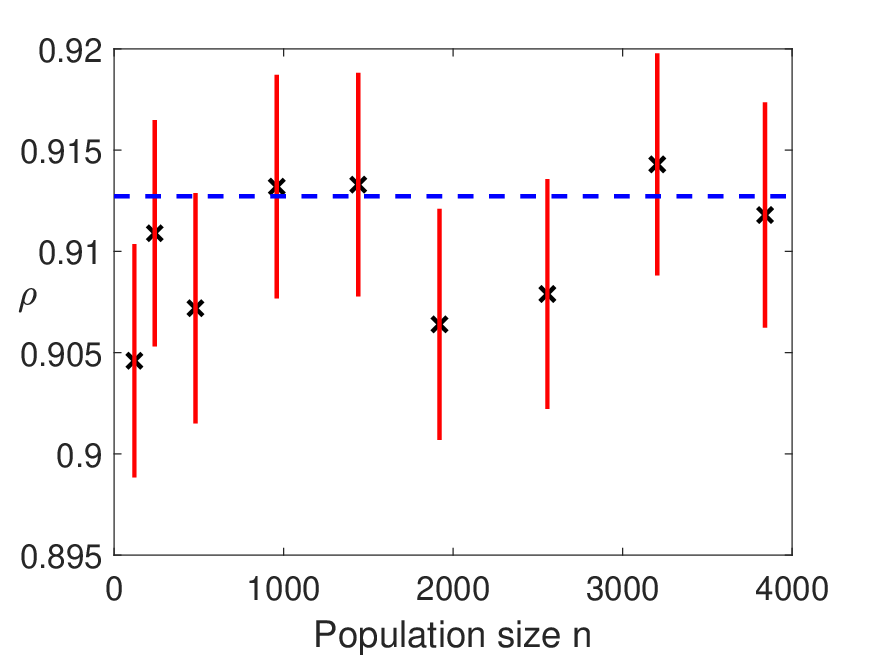} &
 \includegraphics[width=\hfigwidth]{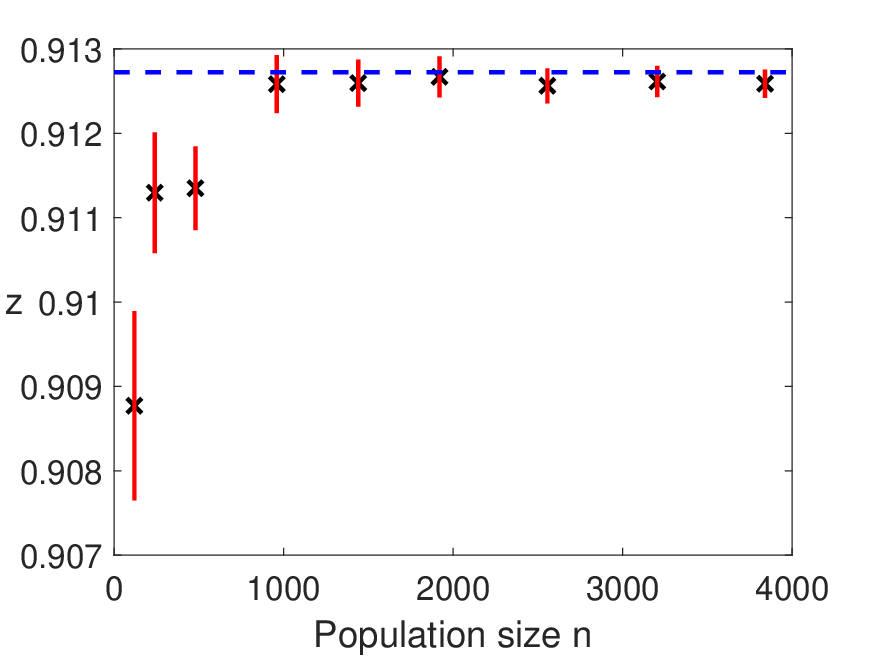}
\end{tabular}
\end{center}
\caption{Comparison of limiting and simulated estimates of the probability of a major outbreak $\rho$ (left column) and the fraction of the population infected by a major outbreak $z$ (right column) based on 10,000 simulations for each choice of $(n, d)$.  The blue dashed horizontal lines depict the asymptotic values.  The black crosses show the estimated values based on simulations, with the red vertical lines giving approximate
$95\%$ confidence intervals.  See text for further details.}
\label{fig:finalsize2}
\end{figure}

In Figure~\ref{fig:finalsize2}, we explore the accuracy of the limiting (as $n \to \infty$) final size $z$ and probability of a major outbreak $\rho$ for finite population sizes.  We assume $\theta=0.4$, $h=4$, $\beta=3$, $\pi_G=0.025$, $\pi_{H|G^c}=0.5$ and $I \equiv 1$, and consider epidemics with $d=1,2$ or $3$ (so $w=4, 8$ or $12$) for various population sizes $n$.  For $d=1$ or $2$, we consider epidemics in populations of size $n=$ 120, 240, 480, 960, 1,440, 1,920, 2,560, 3,200 and 3,840.  The population size $n$ needs to be a multiple of $w=hd$, so when $d=3$, we replaced 2,560 and 3,200 by 2,556 and 3,204, respectively.  The cut-offs for distinguishing between minor and major epidemics were determined by inspection of histograms of simulated final sizes.  In all cases with $n > 480$ a cut-off of $200$ was used.  For $d=1$, the cut-offs for $n=120, 240, 480$ were $36, 75, 150$; for $d=2$, they were $50, 100, 150$; and for $d=3$, they were $25, 50, 100$.  For each fixed $(n,d)$, $n_{\rm sim}=10,000$ epidemics were simulated and an estimate $\hat{\rho}$ of the probability of a major outbreak obtained by the fraction of simulations whose final size was $\ge$ the cut-off, together with an approximate $95\%$ confidence interval for $\rho$ given by $\hat{\rho} \pm 1.96 \sqrt{\hat{\rho}(1-\hat{\rho})/n_{\rm sim}}$.  Independently, 10,000 major outbreaks were simulated, using the given cut-off.  Let $\hat{z}$ and $\hat{\sigma}^2$ be respectively the sample mean and variance of the fraction of the population infected in these simulated epidemics.  Then
$z$ is estimated by $\hat{z}$, with an approximate $95\%$ confidence interval given by $\hat{z} \pm 1.96 \hat{\sigma}/\sqrt{n_{\rm sim}}$.  Note that $\rho=z$ since the infectious period $I$ is non-random.

It is clear from Figure~\ref{fig:finalsize2} that the limiting $\rho$ and $z$ provide good approximations even for moderately sized $n$.  The approximation for the final size $z$ is better than that for the major outbreak probability $\rho$ (note the different scales on the $y$-axes).  In all cases, the limiting $z$ overestimates the final size for finite $n$.  That is not the case for major outbreak probability $\rho$ but note the wider confidence intervals.

\begin{figure}
\begin{center}
\begin{tabular}{ccc}
 \includegraphics[width=\hfigwidth]{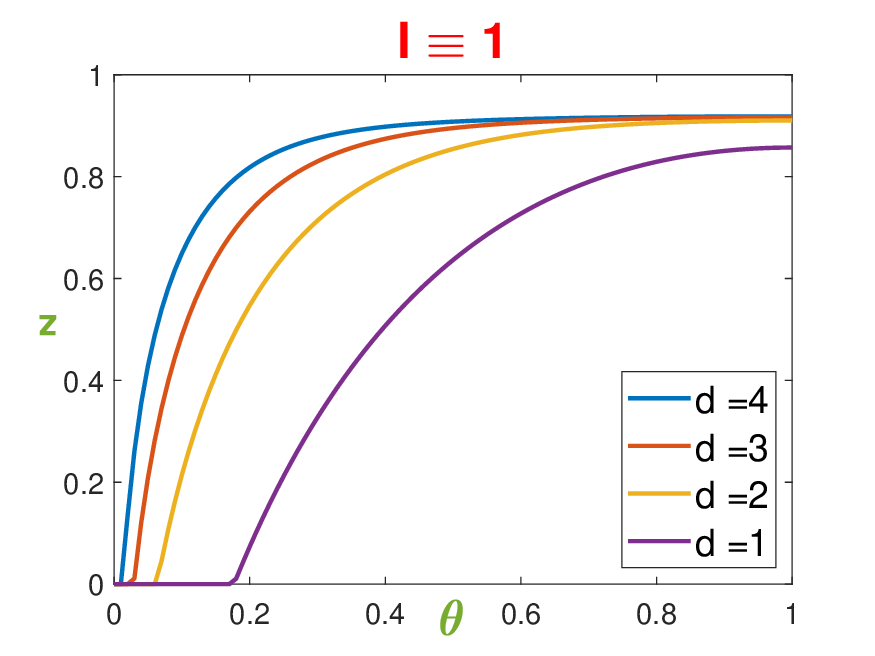} &
 \includegraphics[width=\hfigwidth]{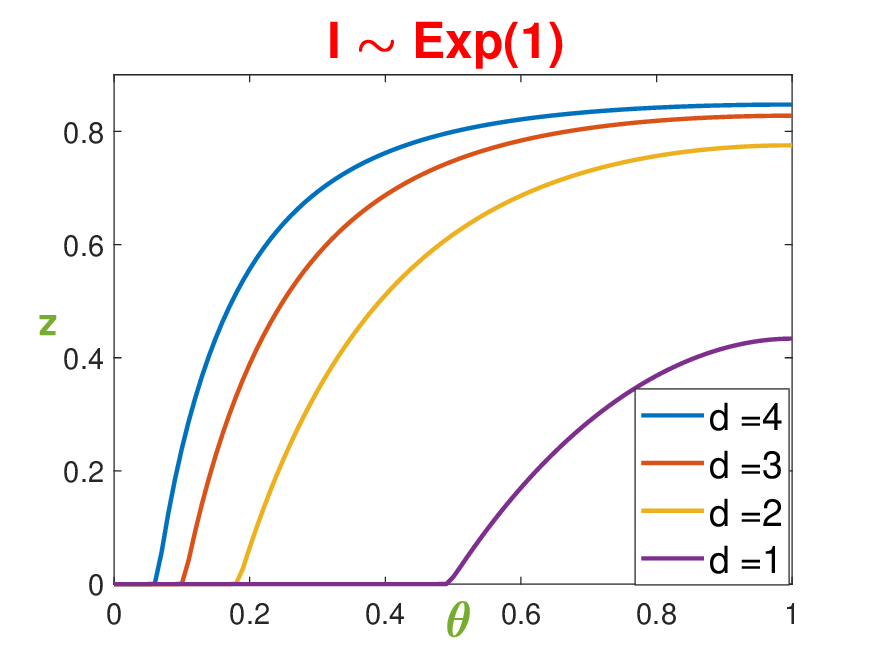}
\end{tabular}
\end{center}
\caption{Dependence of the limiting final size, $z$, on the probability an individual is a mover, $\theta$, for fixed $d=1,2,3,4$ when $h=3$, $\beta=3$, $\pi_G=0.025$ and $\pi_{H|G^c}=0.5$; in the left panel $I \equiv 1$ and in the right panel $I \sim {\rm Exp}(1)$.}
\label{fig:finalsize3}
\end{figure}

In Figure~\ref{fig:finalsize3}, we investigate the dependence of the dependence of the limiting final size, $z$, on the probability an individual is a mover, $\theta$, the workplace to household size ratio, $d$, and whether the infectious period  is non-random or exponential.  We assume $h=3$, $\beta=3$, $\pi_G=0.025$ and $\pi_{H|G^c}=0.5$.  In the left panel $I \equiv 1$ and in the right panel $I \sim {\rm Exp}(1)$.  First note that $z$ is increasing in both $\theta$ and $d$.  Secondly, for fixed $d$, $z$ is fairly constant with $\theta$ when $\theta$ is close to one.  This is particularly the case when the infectious period $I$ is non-random and $d \ge 2$, when $z$ is close to being constant for $\theta \in [0.5,1]$.  Recall that $\theta=1$ yields the standard households-workplaces model analysed for example in Ball and Neal~\cite{BN02} and Pellis et al.~\cite{Pellis09}, in which the population is partitioned independently into households and into workplaces.  Thus, under these circumstances, the standard model in which there is no overlap between the two structures provides a good approximation to the more realistic model having overlap.  Thirdly, note that for the parameter values in Figure~\ref{fig:finalsize3}, $z=0$ when $\theta=0$.  If $\theta=0$, every workplace consists of $d$ complete households, so it is not possible for the disease to spread outside the workplace of the initial case via local infection.  As $\theta$ increases from 0, with all other parameters held constant, there is a critical value of $\theta$, $\hat{\theta}_d$ say, when $R_*=1$ and the epidemic is critical.  Observe that $\hat{\theta}_d$ decreases with $d$ and in an obvious notation, $z'(\hat{\theta}_d)$ is large unless $d$ is small.  Fourthly, given a fixed mean, the choice of the infectious period distribution can have an appreciable impact on epidemic properties.  Note that $z$ is greater when $I \equiv 1$ than in the corresponding model with $I \sim {\rm Exp}(1)$.  More strikingly, $\hat{\theta}_d$ can be a lot smaller for the model with $I \equiv 1$, the difference being particularly noticeable when $d=1$.

\section{Derivations and proofs when $\pi_G>0$}
\label{sec-derivations}
\subsection{Introduction}
Let $\En$ denote the epidemic model defined in Section~\ref{sec-model}.  When $\pi_G>0$, the epidemic model $\En$ is a special case of the general framework of Ball and Neal~\cite{BN02} for SIR epidemics in populations with two levels of mixing.  We first outline that general framework in Section~\ref{sec-genframework}, introducing key concepts of \emph{local infectious clumps}, \emph{local susceptibility sets} and the \emph{severity} of a local infectious clump, and explain heuristically how the limiting severity of a local infectious clump as $n \to \infty$ is related to the event of a major outbreak and the limiting size of a local susceptibility set is related to the final outcome of a major outbreak.  In the present model, the structure of local infectious clumps and local susceptibility sets is far more complicated than in the examples considered in Ball and Neal~\cite{BN02}, and consequently proofs and calculations of limiting major outbreak probability and final size are considerably more difficult.   Local infectious clumps and local susceptibility sets for $\En$ are analysed in Section~\ref{subsec-susc}.  Crucial to this analysis are \emph{complexes}, defined in Section~\ref{subsubsec-complex}, which are associated with workplaces but contain also households of movers who moved out of that workplace.  The analysis of the limiting local infectious clump as $n \to \infty$ is simpler when $I \equiv 1$.  This is treated in Section~\ref{subsubsec-clump}, with the extension to the case of general $I$ and calculation of the limiting major outbreak probability $\rho$ being given in Section~\ref{subsec-outbreakprob}. Local susceptibility sets are analysed in Section~\ref{subsubsec-susset}, where calculations of the threshold parameter $R_*$ and the limiting final size of a major outbreak $z$ are described.  Numerical implementation of the calculation of $R_*$, $z$ and $\rho$ is considered briefly in Section~\ref{subsec-numcomp}.  Finally, proofs of Theorems~\ref{prop-pmajor} and~\ref{prop-finsize} for $\pi_G>0$ are given in Sections~\ref{subsec-majorprobproof} and~\ref{subsec-finsizeproof}, respectively.

\subsection{General framework}
\label{sec-genframework}
We derive the properties of the model by using the general framework of Ball and Neal~\cite{BN02} for SIR epidemics in populations with two levels of mixing. We first outline the key points of that framework in the current setting.  Label the individuals in the population $1, 2, \dots, n$ and let $\mathcal{N}=\{1,2,\dots, n\}$.  (The precise labelling is not important at present.)  We construct on $\mathcal{N}$ a random directed graph $\mathcal{G}^{(n)}$ of potential \emph{local} (i.e.~household or workplace) infectious contacts as follows.  For each $(i,j) \in \mathcal{N}^2$, with $i \ne j$, there is a directed edge from $i$ to $j$ in $\mathcal{G}^{(n)}$
if and only if $i$, if infected, makes infectious contact with $j$.  More precisely, let $I_i$ denote the length of $i$'s infectious period.  Then, given $I_i$, there is a directed edge from $i$ to each individual in $i$'s household independently with probability $1-\exp(-\beta_H' I_i)$, where $\beta_H'=\beta_H/(h-1)$, and a directed edge from $i$ to each individual in $i$'s workplace independently with probability $1-\exp(-\beta_W' I_i)$, where $\beta_W'=\beta_W/(w-1)$.  If an individual belongs to both $i$'s household and workplace, then the directed edge is present with probability $1-\exp(-(\beta_H'+\beta_W')I_i)$.

We use $\mathcal{G}^{(n)}$ to define for each $i \in \mathcal{N}$ their \emph{local infectious clump} and \emph{local susceptibility set} as follows.  For $(i,j) \in \mathcal{N}^2$, write $i \leadsto j$ if and only if there is a chain of directed edges from $i$ to $j$, with the convention that $i \leadsto i$.  For $i \in \mathcal{N}$, let $\Cnmc_i=\{j \in \mathcal{N}:i \leadsto j\}$ and $\Snmc_i=\{j \in \mathcal{N}: j \leadsto i\}$ be $i$'s local infectious clump and local susceptibility set, respectively.  Let $\mathcal{E}^{(n)}$ denote the epidemic model described in Section~\ref{sec-epidemicmodel}.  A realisation of the final outcome of $\mathcal{E}^{(n)}$ can be constructed using $I_1, I_2, \dots, I_n$ and $\mathcal{G}^{(n)}$ by attaching to each $i \in \mathcal{N}$ an independent Poisson process having rate $\beta_G$ giving the times individuals make global contacts if they are infectious. (As in Section~\ref{sec-epidemicmodel}, each global contact is with an individual chosen independently and uniformly from $\mathcal{N}$.) Observe that if $i$ is infected by a global contact in $\mathcal{E}^{(n)}$, then each member of $\Cnmc_i$ is necessarily infected in $\mathcal{E}^{(n)}$.  Similarly, if $i$ is an initial susceptible in $\mathcal{E}^{(n)}$, then $i$ is infected in $\mathcal{E}^{(n)}$ if and only if a member of $\Snmc_i$ is infected by a global contact or $\Snmc_i$ contains the initial infective.

For $i \in \mathcal{N}$, let $\Cn_i=|\Cnmc_i|$ and $\Sn_i=|\Snmc_i|$ be the sizes of $i$'s local infectious clump and local susceptibility set, respectively, and let $\An_i=\sum_{j \in \Cnmc_i} I_j$ be the severity of $i$'s local infectious clump.  Note that $\Cn_1, \Cn_2, \dots \Cn_n$ are identically distributed, with distribution that depends on $n$ (unless $\theta=0$), but they are not independent.  A similar comment holds for $\Sn_1, \Sn_2, \dots, \Sn_n$ and $\An_1, \An_2, \dots, \An_n$.  Let $C^{(n)}, S^{(n)}$ and $A^{(n)}$ be random variables distributed as $\Cn_1, \Sn_1$ and $\An_1$, respectively.  We show below that there exists random variables $C, S$ and $A$ such that
\begin{equation}
\label{equ:CSAconv}
C^{(n)} \convD C, \quad S^{(n)} \convD S \quad \text{and} \quad A^{(n)} \convD A \quad \text{as } n \to \infty.
\end{equation}
The distributions of $C, S$ and $A$ are defective, with a mass at infinity, if the epidemic $\mathcal{E}^{(n)}$ in the absence of global infection (i.e.~with $\beta_G=0$) is supercritical in the limit $n \to \infty$.

Suppose that $\P(C < \infty)=1$ and recall that the epidemic $\mathcal{E}^{(n)}$ is initiated by a single individual being externally infected.  Then, for large $n$, the early stages of $\mathcal{E}^{(n)}$ can be approximated by a branching process of local infectious clumps, since the probability that clumps arising from distinct global contacts have a non-empty intersection is small.   Let $R^{(n)}$ be the number of global contacts that emanate from a typical local infectious clump.  Then $R^{(n)}\sim{\rm Po}(\beta_G A^{(n)})$, i.e.~$R^{(n)}$ is distributed as a mixed Poisson distribution with random mean $\beta_G A^{(n)}$. Thus, $\E[R^{(n)}]=\beta_G \E[A^{(n)}]$.  Moreover, since $I_j$ is independent of the event $\{j \in \Cnmc_i\}$ and $\E[I_j]=1$,
\begin{equation*}
\E[\An_i]=\E\left[\sum_{j=1}^n I_j 1_{\{j \in \Cnmc_i\}}\right] = \sum_{i=1}^n \E[I_j] \P(j \in \Cnmc_i)=\sum_{j=1}^n \E\left[1_{\{j \in \Cnmc_i\}}\right]=\E[\Cn_i].
\end{equation*}
Hence, $\E[R^{(n)}]=\beta_G \E[C^{(n)}]$.  Now $R^{(n)} \convD R$ as $n \to \infty$, where $R \sim {\rm Po}(\beta_G A)$, so the process of local infectious clumps in $\mathcal{E}^{(n)}$ can be approximated by a
Galton-Watson branching process $\mathcal{B}$, having a single ancestor and offspring distribution $R$.

Let $\xi$ be the extinction probability of $\mathcal{B}$.  Then $\xi$ is given by the smallest solution in $[0,1]$ of $f_R(s)=s$. 
Conditioning on $A$,
\begin{equation}
\label{equ:Rpgf}
f_R(s)=\E[\E[s^R|A]]=\E[\exp(-\beta_G A(1-s))]=\phi_A(\beta_G(1-s)).
\end{equation}
Further, $\xi<1$ if and only if $R_*>1$, where
\begin{equation}
\label{equ:Rstar}
R_* = \E[R] = \beta_G \E[A]= \beta_G \E[C].
\end{equation}
Thus, $R_*$ serves as a threshold parameter for the epidemic $\mathcal{E}^{(n)}$, in that as $n \to \infty$ a major outbreak occurs with non-zero probability if and only if $R_*>1$.  Note that $\rho$ in Theorem~\ref{prop-pmajor} is given by $1-\xi$.

The above assumes that $\P(C < \infty)=1$.  If $\P(C < \infty)<1$, then the branching process $\mathcal{B}$ is defined as above but now $R_*=\infty$ and an individual may have infinitely many offspring.  The
branching process $\mathcal{B}$ can still be used to approximate the early stages of $\mathcal{E}^{(n)}$ for large $n$, since an individual having infinitely many offspring in $\mathcal{B}$ necessarily corresponds to a major outbreak in $\mathcal{E}^{(n)}$.  The above heuristic argument is made rigorous in the proof of Theorem~\ref{prop-pmajor} in Section~\ref{subsec-majorprobproof}.

Suppose that $n$ is large and a major outbreak occurs.  Let $z$ be the fraction of the population that are infected.  Since $\E[I]=1$, the total severity of the epidemic is approximately $nz$, so the probability a given individual avoids global contact throughout the epidemic is approximately
\begin{equation}
\label{equ:pi}
\pi=\exp\left(-\frac{\beta_G}{n} nz\right)=\exp(-\beta_G z).
\end{equation}
The probability an individual chosen uniformly at random from the population is infected by the epidemic is given by $z$.  Recall that such an individual, if initially susceptible, avoids infection from the epidemic if and only if (i) all members of its susceptibility set avoid global infection and (ii) its susceptibility set does not contain the initial infective.  The former event has asymptotic probability $0$ if $S=\infty$ and the latter event has asymptotic probability $1$ if $S<\infty$.  For large $n$, individuals avoid global infection approximately independently, so $1-z \approx f_S(\pi)=\E[\pi^S]$.  In Section~\ref{subsec-finsizeproof}, we show that these approximations are exact in the limit as $n \to \infty$ so, using~\eqref{equ:pi}, $z$ satisfies
\begin{equation}
\label{equ:z}
1-z=f_S({\rm e}^{-\beta_G z}).
\end{equation}

Note that $s=1-z$ satisfies the equation governing the extinction probability of a Galton-Watson process having offspring distribution ${\rm Po}(\beta_G S)$.  Further, by exchangeability
\begin{equation}
\label{equ:ECnequalsESn}
\E[C^{(n)}] =\frac{1}{n}\E\left[\sum_{i=1}^n \sum_{j=1}^n 1_{\{i \leadsto j\}}\right]=\frac{1}{n}\E\left[\sum_{j=1}^n \sum_{i=1}^n 1_{\{j \leadsto i\}}\right]= \E[S^{(n)}],
\end{equation}
so letting $n \to \infty$ yields $\beta_G \E[S]=\beta_G \E[C]=R_*$.  It follows from standard branching process theory that if $R_* \le 1$, then $z=0$ is the only solution of~\eqref{equ:z}, while if $R_*>1$, there is a unique second solution in $(0,1]$.  When $R_*>1$, this second solution gives the fraction of the population that is infected by a major outbreak in the limit $n \to \infty$.

In Section~\ref{subsubsec-susset}, we show that $S$, the size of a typical local susceptibility set in the limit $n \to \infty$, is distributed as the total progeny of a $3$-type branching process.  The PGF of $S$ is given in~\eqref{equ:progenypgfS} and $R_*$ is given in~\eqref{equ:RstarS}.
When $I \equiv 1$, the possible directed edges in $\mathcal{G}^{(n)}$ are present independently, with probability $1-\exp(-\beta_H')$ if the two individuals are in the same household but not the same workplace, $1-\exp(-\beta_W')$ if they are in the same workplace but not the same household, and
$1-\exp(-(\beta_H'+\beta_W')$ if they are in both the same household and the same workplace.  It follows that $\Cn \eqd \Sn$, so $C \eqd S$.  Note that in this case $\An=\Cn$, whence $A \eqd C \eqd S$, and~\eqref{equ:Rpgf} and~\eqref{equ:z} imply that $\rho=z$.  When $I$ is non-constant,  $C$ is distributed as the total progeny of a $3(h+w-2)$-type branching process and $A$ as the sum of the lifetimes of all individuals in that branching process.  The derivation of the Laplace transform $\phi_A(\nu)$ in that case is outlined in Section~\ref{subsec-outbreakprob}.

Although this section is concerned with the case $\pi_G>0$, we note here that when $\pi_G=0$, we have $\rho=\P(C=\infty)$ and $z=\P(S=\infty)$.

\subsection{Local infectious clumps and susceptibility sets}\label{subsec-susc}
\subsubsection{Complexes}
\label{Complex}
\label{subsubsec-complex}

In order to apply the theory outlined in Section~\ref{sec-genframework}, we need to obtain the distribution of the random variables $C, A$ and $S$, describing the size and severity of a typical local infectious clump and the size of a typical local susceptibility set, respectively, in the limit as $n \to \infty$.  For this purpose, we first define a {\it complex}.  There are $m = n/w$ complexes in the network, where $m$ is the number of workplaces. Each complex is derived from a workplace containing $d$ households with individuals belonging to 1 or 2 complexes depending on whether they are a remainer (belong to a single complex) or a mover (belongs to two complexes). We label the workplaces $1,2,\dots, m$ and associate the $i^{th}$ complex with the $i^{th}$ workplace. There are 3 types of individuals in the $i^{th}$ complex:
\begin{itemize}
\item Movers who start in workplace $i$  (their household belongs to workplace $i$) but move to a different workplace.
\item  Remainers who start in the workplace $i$ (their household belongs to workplace $i$)  and remain in that workplace.
\item Movers who start in a different workplace (their household is part of a different workplace) and move into workplace $i$.
\end{itemize}

We consider complex $i$. We label the households in complex $i$, $(i,1), (i,2), \ldots, (i,d)$. Then complex $i$ is comprised of $2 d+1$ groups of individuals. These are:
\begin{itemize}
\item Group $2j-1$ $(j=1,2,\ldots,d)$. These are the individuals in household $(i,j)$ who remain in workplace $i$.
\item Group $2j$ $(j=1,2,\ldots,d)$. These are the individuals in household $(i,j)$ who move to a different workplace.
\item Group $2d+1$. The individuals who move into workplace $i$.
\end{itemize}
Therefore in complex $i$, household-based infections take place in household $(i,j)$ amongst the individuals in groups $2j-1$ and $2j$ and workplace-based infections involve groups $j=1,3,\ldots, 2d+1$ (the odd numbered groups).  We note that group $2j-1$ $(j=1,2,\ldots,d)$ can be empty if all members of household $(i,j)$ move to different workplaces. We still associate the household $(i,j)$ with workplace $i$.

In Figure \ref{fig:complex} we illustrate how complexes and the population are constructed using an example where $h=d=2$ and $m=4$ giving $w=4$ and $n =16$.

When $\theta=0$, everyone is a remainer, the complexes do not intersect and each complex consists of a workplace that is partitioned into $d$ households.  The distributions of $C^{(n)}, A^{(n)}$ and $S^{(n)}$ are each independent of $n$, so $C^{(n)} \eqd C, A^{(n)} \eqd A$ and $S^{(n)} \eqd S$ for all $n$.  The threshold parameter $R_*$ and the Laplace transform $\phi_A(\nu)$ can be computed, at least in principle, using exact results for the final outcome of stochastic multitype SIR epidemics given in Ball~\cite{Ball86}, Section 3, or Picard and Lef\`evre~\cite{PL90}, Section 4.  The probability mass function of $S$ can be computed in principle  using Ball~\cite{Ball19}, Lemma 4.1.  In practice, these computations are possible only for suitably small $d$.
In the following, we assume that $\theta>0$.

 \begin{figure}
(a)
\begin{center}
\begin{tikzpicture}

\draw[red] (0,0)  rectangle (2.5,2.5);

\draw[green] (0.75,1.25) ellipse (0.4 and 0.8);
\draw[green] (1.75,1.25) ellipse (0.4 and 0.8);

\draw[red] (3,0)  rectangle (5.5,2.5);

\draw[green] (3.75,1.25) ellipse (0.4 and 0.8);
\draw[green] (4.75,1.25) ellipse (0.4 and 0.8);

\draw[red] (6,0)  rectangle (8.5,2.5);

\draw[green] (6.75,1.25) ellipse (0.4 and 0.8);
\draw[green] (7.75,1.25) ellipse (0.4 and 0.8);

\draw[red] (9,0)  rectangle (11.5,2.5);

\draw[green] (9.75,1.25) ellipse (0.4 and 0.8);
\draw[green] (10.75,1.25) ellipse (0.4 and 0.8);

\node[black] at (1.25,-0.5) {A};
\node[black] at (4.25,-0.5) {B};
\node[black] at (7.25,-0.5) {C};
\node[black] at (10.25,-0.5) {D};

\node[circle, black,fill] at (0.75,0.9) {\tiny 0};
\node[circle, black,fill] at (0.75,1.6) {\tiny 0};

\node[circle, black,draw] at (1.75,0.9) {\tiny 1};
\node[circle, black,draw] at (1.75,1.6) {\tiny 2};

\node[circle, black,fill] at (3.75,0.9) {\tiny 0};
\node[circle, black,draw] at (3.75,1.6) {\tiny 3};

\node[circle, black,fill] at (4.75,0.9) {\tiny 0};
\node[circle, black,draw] at (4.75,1.6) {\tiny 4};

\node[circle, black,fill] at (6.75,0.9) {\tiny 0};
\node[circle, black,fill] at (6.75,1.6) {\tiny 0};

\node[circle, black,fill] at (7.75,0.9) {\tiny 0};
\node[circle, black,draw] at (7.75,1.6) {\tiny 5};

\node[circle, black,fill] at (9.75,0.9) {\tiny 0};
\node[circle, black,draw] at (9.75,1.6) {\tiny 6};

\node[circle, black,draw] at (10.75,0.9) {\tiny 7};
\node[circle, black,draw] at (10.75,1.6) {\tiny 8};

\end{tikzpicture}

\end{center}

\vspace*{0.1cm}

(b)
\begin{center}
\begin{tikzpicture}

\draw[red] (0,0)  rectangle (2.5,2.5);

\draw[green] (0.75,1.25) ellipse (0.4 and 0.8);
\draw[green] (3.25,1.25) ellipse (0.4 and 0.8);

\node[black] at (1.25,-0.5) {A};

\node[circle, black,fill] at (0.75,0.9) {\tiny 0};
\node[circle, black,fill] at (0.75,1.6) {\tiny 0};

\node[black,draw] at (1.75,0.9) {\tiny 4};
\node[black,draw] at (1.75,1.6) {\tiny 6};

\node[circle, black,draw] at (3.25,0.9) {\tiny 1};
\node[circle, black,draw] at (3.25,1.6) {\tiny 2};

\draw[red] (6,0)  rectangle (8.5,2.5);

\draw[green] (8.5,0.7) ellipse (0.8 and 0.4);
\draw[green] (8.5,1.8) ellipse (0.8 and 0.4);

\node[black] at (7.25,-0.5) {B};

\node[circle, black,fill] at (8.2,0.7) {\tiny 0};
\node[circle, black,fill] at (8.2,1.8) {\tiny 0};

\node[black,draw] at (6.75,0.7) {\tiny 2};
\node[black,draw] at (6.75,1.8) {\tiny 7};

\node[circle, black,draw] at (8.8,0.7) {\tiny 3};
\node[circle, black,draw] at (8.8,1.8) {\tiny 4};

\draw[red] (0,-4)  rectangle (2.5,-1.5);

\draw[green] (0.75,-2.75) ellipse (0.4 and 0.8);
\draw[green] (2.5,-3.3) ellipse (0.8 and 0.4);

\node[black] at (1.25,-4.5) {C};

\node[circle, black,fill] at (0.75,-3.1) {\tiny 0};
\node[circle, black,fill] at (0.75,-2.4) {\tiny 0};

\node[circle, black,fill] at (2.2,-3.3) {\tiny 0};
\node[circle, black,draw] at (2.8,-3.3) {\tiny 5};

\node[black,draw] at (2,-2.3) {\tiny 8};

\draw[red] (6,-4)  rectangle (8.5,-1.5);

\draw[green] (8.5,-3.3) ellipse (0.8 and 0.4);

\draw[green] (10,-2.75) ellipse (0.4 and 0.8);

\node[black] at (7.25,-4.5) {D};

\node[circle, black,fill] at (0.75,-3.1) {\tiny 0};
\node[circle, black,fill] at (0.75,-2.4) {\tiny 0};

\node[circle, black,fill] at (8.2,-3.3) {\tiny 0};
\node[circle, black,draw] at (8.8,-3.3) {\tiny 6};

\node[circle, black,draw] at (10,-2.4) {\tiny 7};
\node[circle, black,draw] at (10,-3.1) {\tiny 8};

\node[black,draw] at (6.75,-2.3) {\tiny 1};
\node[black,draw] at (7.75,-2.3) {\tiny 3};
\node[black,draw] at (6.75,-3.1) {\tiny 5};

\end{tikzpicture}
\end{center}

\vspace*{0.1cm}

(c)
\begin{center}
\begin{tikzpicture}

\draw[red] (0,0)  rectangle (2.5,2.5);

\draw[red] (3,1.75)  rectangle (5.5,4.25);

\draw[red] (3,-1.25)  rectangle (5.5,1.25);

\draw[red] (6,0)  rectangle (8.5,2.5);

\node[black] at (-0.5,1.25) {A};
\node[black] at (2.5,4) {B};
\node[black] at (9,1.25) {C};
\node[black] at (2.5,-1) {D};

\draw[green] (0.75,1.25) ellipse (0.4 and 0.8);
\draw[rotate around={135:(2.75,2.75)},green] (2.75,2.75) ellipse (0.65 and 1.3);

\draw[rotate around={45:(2.75,0)},green] (2.75,0) ellipse (0.6 and 1.2);

\draw[rotate around={45:(5.75,2.75)},green] (5.75,2.75) ellipse (0.65 and 1.3);

\draw[rotate around={135:(5.75,0)},green] (5.75,0) ellipse (0.6 and 1.2);
\draw[green] (7.75,1.25) ellipse (0.4 and 0.8);

\draw[green] (3.75,1.5) ellipse (0.475 and 0.95);
\draw[green] (4.75,1.5) ellipse (0.475 and 0.95);

\node[circle, black,fill] at (0.75,0.9) {\tiny 0};
\node[circle, black,fill] at (0.75,1.6) {\tiny 0};

\node[circle, black,fill] at (3.3,-0.5) {\tiny 0};
\node[circle, black,fill] at (3.3,3.25) {\tiny 0};

\node[circle, black,draw] at (5.2,-0.5) {\tiny 5};
\node[circle, black,draw] at (5.2,3.25) {\tiny 7};

\node[circle, black,draw] at (2.2,0.5) {\tiny 6};
\node[circle, black,draw] at (2.1,2.1) {\tiny 4};

\node[circle, black,fill] at (6.3,0.5) {\tiny 0};
\node[circle, black,draw] at (6.4,2.1) {\tiny 8};

\node[circle, black,fill] at (7.75,0.9) {\tiny 0};
\node[circle, black,fill] at (7.75,1.6) {\tiny 0};

\node[circle, black,draw] at (3.75,2.1) {\tiny 2};
\node[circle, black,draw] at (3.75,0.9) {\tiny 1};

\node[circle, black,fill] at (4.75,2.1) {\tiny 0};
\node[circle, black,draw] at (4.75,0.9) {\tiny 3};

\end{tikzpicture}
\end{center}

\caption{\footnotesize Illustration of the construction of a population of $n=16$ individuals divided into 4 workplaces, labelled $A$, $B$, $C$, $D$, of size $4(=w)$ with each workplace consisting of $2(=d)$ households of size $2(=h)$. Red squares denote workplaces and green ellipses denote households. \break
(a) The initial workplaces with solid circles representing remainers and open circles representing movers with movers numbered. There are 8 movers and 8 remainers.
\break
(b) The construction of complexes. Open squares represent the places created in a workplace by movers and the number in the square denotes the mover who fills the place.  Note that in complexes A and D there are households consisting completely of movers with the household residing outside the workplace. \break
(c) The final population structure created from the complexes.}
\label{fig:complex}
\end{figure}

\subsubsection{Local infectious clumps}
\label{subsubsec-clump}
Suppose first that the infectious period $I$ is non-random, so $\P(I=1)=1$. We construct the local infectious clump of a typical individual, $i^*$ say, on a generation basis by considering local epidemic spread within successive complexes.  Note that apart from the initial generation, a local epidemic {\it enters} a complex through one of two ways: either a mover in one of groups $2, 4, \ldots, 2d$ is infected in their workplace or a mover in group $2d+1$ is infected in their household.  The initial individual $i^*$ can be any member of a complex.

We observe that the groups $1,3,\ldots, 2d-1$ are exchangeable, as are the  groups $2,4,\ldots, 2d$. Therefore we only need to consider the behaviour of the local spread in a complex in the cases where the initial infective in a complex is in group 1, 2 or $2d+1$.  Let $R$ denote that the initial infective is in group $1$ (or $3,5,\ldots,2d-1$), that is, the individual is a remainer. Let $H$ denote that the initial infective is in group $2$ (or $4,6,\ldots,2d$), the individual belongs to complex $i$ through their household. Let $W$ denote that the initial infective is in group $2d+1$,  the individual belongs to complex $i$ through the workplace.

The initial infective in a complex being of type $X$ is informative about the structure of the complex. For example, if $X=W$, then we know that there is at least one mover in the workplace of the complex. Without loss of generality, if the initial individual is of type $R$ or $H$ we assign them to group 1 or 2, respectively. Let $M_1 \sim {\rm Bin} (n-1,\theta)$ and for $j=2,3,\ldots,d$, let $M_j \sim {\rm Bin} (n, \theta)$. The number of individuals in groups $2j-1$ and $2j$ $(j=2,3,\ldots,d)$ are $n-M_j$ and $M_j$, respectively, regardless of the type of the initial infective. If $X=R$, then the number of individuals in groups 1, 2 and $2d+1$ are $n-M_1$, $M_1$ and $\sum_{l=1}^d M_l$, respectively, otherwise there are  $n-1-M_1$, $1+M_1$ and $1+\sum_{l=1}^d M_l$ individuals in groups 1, 2 and $2d+1$, respectively.

If $n$ is large, the probability that, in the early stages of the construction of an infectious clump, there is a complex that contains movers who either moved into the complex from the same workplace or moved out of their original workplace to the same workplace is small.  Thus, we approximate the size $C^{(n)}$ of the clump by assuming that for each complex, all the members from groups $2,4,\ldots, 2d$ are in different workplaces and all members of group $2d+1$ are in different households.  Further, we assume that apart from the initial infective in a complex, all movers in a complex have not been involved previously in the construction of the clump.  Under these assumptions, $C^{(n)}$ is approximated by the total progeny of a 3-type branching process, denoted by $\mathcal{B}_C$, which we now describe.

Consider the spread of a local epidemic within a typical complex defined as above.  For $X, Y=R, H, W$, let $Z_{XY}$ be defined as follows. Let $Z_{XR}$ denote the number of remainer individuals (excluding the initial infective if $X=R$) infected in a complex when there is a single initial infective in group type $X$.
Let $Z_{XH}$ denote the number of individuals in group $2d+1$ (excluding the initial infective if $X=W$) infected in the complex when there is a single infective in group type $X$. Individuals infected in group $2d+1$ will be infected in the workplace and will start an epidemic in a new complex through their household. In other words, in that other complex they will belong to group $2, 4, \ldots, 2d$.
Similarly, let $Z_{XW}$ denote the number of individuals in groups $2, 4, \ldots, 2d$ (excluding the initial infective if $X=H$) infected in the complex when there is a single infective in group type $X$. Individuals infected in group  $2, 4, \ldots, 2d$ will be infected in the household and will start an epidemic in a new complex through their workplace. In other words, in that other complex they will belong to group $2d+1$.

If $i^*$, the initial member of the local infectious clump is a remainer, then the initial individual in $\mathcal{B}_C$ has type $R$ and its offspring is distributed as $(Z_{RR}, Z_{RH}, Z_{RW})$.
If $i^*$ is a mover, then $i^*$ belongs to two complexes and its offspring is distributed as $(Z_{HR}, Z_{HH}, Z_{HW})+(Z_{WR}, Z_{WH}, Z_{WW})$, where the two random vectors are independent since $I$ is non-random.  All subsequent individuals in the branching process will be of type $H$ or $W$, and their offspring are distributed as $(Z_{HR}, Z_{HH}, Z_{HW})$ or $(Z_{WR}, Z_{WH}, Z_{WW})$, respectively.  Note that individuals in the branching process $\mathcal{B}_C$ correspond to complexes in the local infectious clump.

Consider the branching process $\hat{\mathcal{B}}_C$, in which there is a single ancestor, who is either of type $H$ or $W$, and the offspring distributions of type-$H$ and type-$W$ individuals are $(Z_{HR}, Z_{HH}, Z_{HW})$ and $(Z_{WR}, Z_{WH}, Z_{WW})$, respectively.  Thus, type-$R$ individuals may exist in $\hat{\mathcal{B}}_C$ but they have no offspring.  Let $\hat{Z}_H$ be the total progeny of $\hat{\mathcal{B}}_C$ (i.e~the total number of individuals infected and not the total number of complexes infected), counting all types and including the ancestor, given that the ancestor is of type $H$, and define $\hat{Z}_W$ similarly.
Also, for $X=R, H, W$ and $0 \le s_1, s_2, s_3 \le 1$, let
\[
g_X(s_1, s_2, s_3)=\E\left[s_1^{Z_{XR}}s_2^{Z_{XH}}s_3^{Z_{XH}}\right].
\]
By taking expectations with respect to the offspring of the ancestor,
\begin{equation}
\label{equ:progenypgfhat}
f_{\hat{Z}_H}(s)=s g_H(s, f_{\hat{Z}_H}(s), f_{\hat{Z}_W}(s)) \qquad \text{and} \qquad f_{\hat{Z}_W}(s)=s g_W(s, f_{\hat{Z}_H}(s), f_{\hat{Z}_W}(s)).
\end{equation}
For each $s \in [0,1]$, the equations~\eqref{equ:progenypgfhat} determine $f_{\hat{Z}_H}(s)$ and $f_{\hat{Z}_W}(s)$.

Let $C$ be the total progeny of $\mathcal{B}_C$, counting all types and including the ancestor, given that the ancestor is a mover with probability $\theta$, otherwise it is a remainer.
Then, conditioning on the type of the ancestor and the spread within its complex(es),
\begin{equation}
\label{equ:progenypgf}
f_C(s)=(1-\theta)sg_R(s, f_{\hat{Z}_H}(s), f_{\hat{Z}_W}(s))+\theta g_H(s, f_{\hat{Z}_H}(s), f_{\hat{Z}_W}(s))f_{\hat{Z}_W}(s).
\end{equation}
(Note, using~\eqref{equ:progenypgfhat}, that for $s>0$, the second term on the right-hand side of~\eqref{equ:progenypgf} may be expressed as $\theta s^{-1}f_{\hat{Z}_H}(s)f_{\hat{Z}_W}(s)$; we use the form in~\eqref{equ:progenypgf} as it is also valid when $s=0$.)

${\newtheorem{lemma-rfclumpconv}[theorem]{Lemma}}$
\begin{lemma-rfclumpconv} \label{lem:rfclumpconv}
For the model with a constant infectious period,
\begin{equation*}
C^{(n)} \convD C \quad \text{as } n \to \infty.
\end{equation*}
\end{lemma-rfclumpconv}

\begin{proof}
We prove the lemma by coupling realisations of the local infectious clump $\mathcal{C}^{(n)}_{i^*}$ of an individual $i^*$ in $\mathcal{E}^{(n)}$ and the branching process $\mathcal{B}_C$, both of which are constructed on a generation basis, taking complexes in the same order.  We start with a realisation of $\mathcal{C}^{(n)}_{i^*}$ and use the local spread within its complexes to give a realisation of $\mathcal{B}_C$ until in $\mathcal{C}^{(n)}_{i^*}$ we encounter a complex which contains a mover whose original workplace is the same as that of another mover used previously in the construction of $\mathcal{C}^{(n)}_{i^*}$.  Let $D^{(n)}$ be the number of complexes used in $\mathcal{C}^{(n)}_{i^*}$ before we encounter such a complex, where $D^{(n)}=\infty$ if no such complex is encountered.  As soon as such a complex is encountered, the offsprings of all subsequent individuals in $\mathcal{B}_C$ are sampled independently of $\mathcal{C}^{(n)}_{i^*}$ from their appropriate distributions.

Note that for $k=1,2,\dots$, if $C^{(n)}_{i^*}=k$ and $D^{(n)}>k+1$, then
$C^{(n)}_{i^*}=C$, since in the construction of $\mathcal{C}^{(n)}_{i^*}$, each new complex necessarily contains at least one new individual.
Hence, for $k=1,2,\dots$,
\begin{align*}
\label{equ:Clek}
\P(C^{(n)} \le k)&=\P(C^{(n)} \le k, D^{(n)}>k+1)+\P(C^{(n)} \le k, D^{(n)} \le k+1)\nonumber\\
&=\P(C \le k, D^{(n)}>k+1)+\P(C^{(n)} \le k, D^{(n)} \le k+1).
\end{align*}
Also,
\begin{equation*}
\label{equ:Zlek}
\P(C \le k)=\P(C \le k, D^{(n)}>k+1)+\P(C \le k, D^{(n)} \le k+1),
\end{equation*}
It follows that
\begin{equation}
\label{equ:PCminusPZ}
|\P(C^{(n)} \le k)-\P(C \le k)| \le \P(D^{(n)} \le k+1).
\end{equation}

Recall that $\theta>0$.  Let $T^{(n)}$ be the total number of movers in $\mathcal{E}^{(n)}$.  A simple application of Chebyshev's inequality yields $\P(T^{(n)}<\theta n/2) \to 0$ as $n \to \infty$.  The maximum size of a complex is $2w=2h d$, so the number of movers in the first $k+1$ complexes is at most $2w(k+1)$.  Given, $T^{(n)} \ge \theta n/2$ the probability that any two given distinct movers come from the same original workplace is bounded above by $(w-1)/(\frac{1}{2}\theta n-1)$.  Thus,
\[
\P(D^{(n)} \le k+1) \le \P(T^{(n)}<\theta n/2)+\binom{2w(k+1)}{2} \frac{w-1}{\frac{1}{2}\theta n-1} \P(T^{(n)}>\theta n/2) \to 0 \quad \text{as } n \to \infty.
\]
Hence, for $k=1,2,\dots$,
\[
\lim_{n \to \infty} \P(C^{(n)} \le k) = \P(C \le k),
\]
and the lemma follows.
\end{proof}

When the infectious period $I$ is random, the local infectious clump $\mathcal{C}^{(n)}_{i^*}$ can still be constructed on a generation basis and coupled to a limiting process as in Lemma~\ref{lem:rfclumpconv}.  However, the limiting process is not a branching process, since the infectious periods of movers infected in a within-complex epidemic (i.e.~individuals who subsequently start epidemics in new complexes) are not independent of the size of that within-complex epidemic.  The limiting process can be made a branching process by further typing, for example, by also typing infected movers by the number they infect in their household/workplace in the next complex, as we describe in Section~\ref{subsec-outbreakprob}.

\subsubsection{Major outbreak probability}\label{subsec-outbreakprob}
Recall from Section~\ref{sec-genframework} that the major outbreak probability $\rho=1-\xi$, where $\xi$ is the smallest root in $[0, 1]$ of $f_R(s)=s$, with $f_R(s)$ being given by~\eqref{equ:Rpgf}.  Suppose first that $I \equiv 1$.  Then, as noted in Section~\ref{sec-genframework}, we have that $A \eqd C$ and $\rho=z$.  We describe the computation of $z$ (for the case of general $I$) in Section~\ref{subsubsec-susset} below, thus enabling $\rho$ to be calculated.  Alternatively, substituting $A \eqd C$ into~\eqref{equ:Rpgf} yields $f_R(s)=f_C({\rm e}^{-\beta_G(1-s)})$, so $f_R(s)$ can be obtained using~\eqref{equ:progenypgf}.
We outline below how the distribution of the limiting severity $A$ can be obtained for general $I$ in a constructive manner, from which $\rho$ can be obtained using \eqref{equ:Rpgf}.  The key building block is the epidemic within a complex and in particular the joint distribution of the severity of the within-complex epidemic and the number of infectious contacts outside the complex made by movers infected within the complex.

For $X=W,H$ and $l=1,2,\ldots, w-1$, let $A^{X,l}$ be the severity, excluding the initial infective, in a complex, where the initial infective in the complex is a mover who belongs to the complex only through its workplace if $X=W$ (household if $X=H$) and the initial infective infects $l$ individuals within the complex (either in the workplace or household depending on whether $X=W$ or $X=H$). Note that if we know the number of infectious contacts the initial infective makes within a complex then the size of the epidemic in the complex
is conditionally independent of the initial infective's infectious period. In addition, for $Y=W,H$ and $k=1,2,\ldots, w-1$, let $Z^{X,l}_{Y,k}$ denote the number of type $(Y,k)$ movers infected in the complex with an initial infective of type-$(X,l)$. A type-$(Y,k)$ individual is a mover who belongs to another complex through their workplace if $Y=W$ (household if $Y=H$) and infects $k$ individuals within their other complex. We are interested in the joint distribution of $(A^{X,l}, Z^{X,l})$, where $\mathbf{Z}^{X,l} = (Z^{X,l}_{H,1},\ldots Z^{X,l}_{H,h-1}, Z^{X,l}_{W,1}, \ldots Z^{X,l}_{W,w-1})$ is a random vector of length $w+h-2$.

For $X=W,H$ and $l=1,2,\ldots, w-1$, let $\hat{A}^{X,l}$ be the severity (sum of infectious periods), excluding the initial infective, in a {\it branch} of the infectious clump, where the initial infective in a branch is a mover who belongs to the initial complex in a branch only through its workplace if $X=W$ (household if $X=H$) and the initial infective infects $l$ individuals within the complex (either in the workplace or household depending on whether $X=W$ or $X=H$). A branch of an infectious clump starts with a mover individual who begins a complex epidemic and consists of all subsequent complex epidemics emanating from the initial infectious complex in the branch. By the (asymptotic) branching process of the growth of an infectious clump it follows that
\begin{align} \label{eq:outbreak:branch}
\hat{A}^{X,l} \eqd A^{X,l} + \sum_{(Y,k)} \sum_{i=1}^{Z^{X,l}_{Y,k}} \hat{A}^{Y,k,i},
\end{align}
where $\hat{A}^{Y,k,1}, \hat{A}^{Y,k,2}, \ldots$ are independent and identically distributed copies of $\hat{A}^{Y,k}$ and sums are equal to 0 if vacuous.

For a given complex structure, i.e.~sizes of groups $1,2, \dots, 2d+1$, $\mathbf{Z}^{X,l}$ is a vector of final state random variables defined on a multitype SIR epidemic and the joint Laplace transform generating function of $(A^{X,l},\mathbf{Z}^{X,l})$ can be computed using Ball and O'Neill~\cite{Ball-ONeill-1999}, Theorem 5.1, or Ball~\cite{Ball19}, Theorem 4.2.  The unconditional joint Laplace transform generating function of $(A^{X,l}, \mathbf{Z}^{X,l})$ can then be obtained by taking expectations with respect to the complex structure.  For fixed $\nu \ge 0$, a set of $w+h-2$ non-linear equations which determine $(\phi_{\hat{A}^{H,1}}(\nu), \ldots, \phi_{\hat{A}^{H,h-1}}(\nu), \phi_{\hat{A}^{W,1}}(\nu), \ldots, \phi_{\hat{A}^{W,w-1}}(\nu))$ can then be derived using~\eqref{eq:outbreak:branch} enabling, at least in principle, the $\phi_{\hat{A}^{X,l}}(\nu)$ to be computed.

Consider the initial individual in the infectious clump. Let $I$ denote its infectious period, and $Q_H$ and $Q_W$ denote the number of individuals the initial infective infects within their household and workplace, respectively. Then $Q_H | I \sim {\rm Bin} (h-1, 1- \exp(-\beta_H I))$ and $Q_W | I \sim {\rm Bin} (w-1, 1- \exp(-\beta_W I))$, where $Q_H$ and $Q_W$ are conditionally independent given $I$. Let $\tilde{A}_M$ denote the severity of an infectious clump where the initial infective is a mover. Conditional on the triple $(I,Q_H, Q_W)$,
\begin{align} \label{eq:outbreak:mover}
\tilde{A}_M \eqd I + \hat{A}^{H,Q_H} + \hat{A}^{W,Q_W},
\end{align}
with $I$, $\hat{A}^{H,Q_H}$ and $\hat{A}^{W,Q_W}$ conditionally independent given $Q_H$ and $Q_W$. (The initial infective is responsible for instigating two branches, originating from each of the complexes to which it belongs.)
By conditioning first on $I$ and then on $(Q_H, Q_W)$, the Laplace transform $\phi_{\tilde{A}_M}(\nu)$ can be computed using~\eqref{eq:outbreak:mover}.

Consider now the case when the initial infective in an infectious clump is a remainer.  Then the initial infective belongs to only one complex; in all subsequent complexes in the clump, the initial infective in the complex is a mover, but in this first complex they are a remainer.  Let $A^{R,(j,l)}$ denote the severity, excluding the initial infective, in a complex, where the initial infective in the complex is a remainer who makes contact with $j$ of its housemates and $l$ of its workplace colleagues. Note that some individuals can belong to both the same household and workplace as the initial infective so that the number of individuals infected by the initial infective can be less than $j+l$. For $Y=W,H$ and $k=1,2,\ldots, w-1$, let $Z^{R,(j,l)}_{Y,k}$ denote the number of type $(Y,k)$ movers infected in the complex with an initial infective of type $(R,j,l)$.
Let $\tilde{A}_R$ denote the severity of an infectious clump where the initial infective is a remainer. Conditional on the triple $(I,Q_H,Q_W)$,
\begin{align} \label{eq:outbreak:remainer}
\tilde{A}_R \eqd I + A^{R,(Q_H,Q_W)} + \sum_{(Y,k)} \sum_{i=1}^{Z^{R,(Q_H,Q_W)}_{Y,k}} \hat{A}^{Y,k,i},
\end{align}
where $ \hat{A}^{Y,k,1},  \hat{A}^{Y,k,2}, \ldots$ are independent and identically distributed copies of $ \hat{A}^{Y,k}$ and are independent of $(I,A^{R,(Q_H,Q_W)})$ given $Z^{R,(Q_H,Q_W)}$.
For each $(j,l)$, the joint Laplace transform generating function of $(A^{R,(j,l)},\mathbf{Z}^{R,(j,l)})$ can be computed similarly to that of $(A^{X,l}, \mathbf{Z}^{X,l})$ above, whence
$\phi_{\tilde{A}_R}(\nu)$ can be computed using~\eqref{eq:outbreak:remainer}.

Finally we note that $A$ is a mixture of $\tilde{A}_M$ and $\tilde{A}_R$ with mixing probabilities $\theta$ and $1-\theta$, so
\begin{equation}
\label{equ:phiA}
\phi_A(\nu)=\theta \phi_{\tilde{A}_M}(\nu)+(1-\theta)\phi_{\tilde{A}_R}(\nu),
\end{equation}
and using~\eqref{equ:Rpgf}, $\rho=1-\xi$, where $\xi$ is the smallest root of $\phi_A(\beta_G(1-s))=s$ in $[0,1]$.

\subsubsection{Local susceptibility sets}
\label{subsubsec-susset}
The local susceptibility set $\mathcal{S}^{(n)}_{i^*}$ of an individual $i^*$ in $\mathcal{E}^{(n)}$ can be constructed on a generation basis in an analogous fashion to the construction of $\mathcal{C}^{(n)}_{i^*}$ in Section~\ref{subsubsec-clump}.
We first construct the susceptibility set of $i^*$ within the one or two complexes to which $i^*$ belongs.  We then construct the susceptibility sets of all movers (apart from $i^*$ if it is a mover) in that susceptibility set in their other complex, and so on.  An equivalent result to Lemma~\ref{lem:rfclumpconv} holds.  Moreover, the limiting process is necessarily a branching process, even if $I$ is random, since when considering the susceptibility $\mathcal{S}_{C,j^*}$ say, of a single individual $j^*$ within a complex $C$, the distribution of $\mathcal{S}_{C,j^*}$ does not depend on $I_{j^*}$.

A similar argument to the proof of Lemma~\ref{lem:rfclumpconv} shows that $\Sn \convD S$ as $n \to \infty$, where $S$ is the total progeny  of a 3-type branching process, denoted by $\mcB_S$, which is defined analogously to $\mcB_C$ in Section~\ref{subsubsec-clump}. For $X, Y=R, H, W$, let $Z_{XY}^S$ be defined similarly to $Z_{XY}$ in Section~\ref{subsubsec-clump} but for the susceptibility set of an individual within a typical complex.  Thus, in the above notation, $Z_{XY}^S$ is the number of type-$Y$ individuals in $\mathcal{S}_{C,j^*}\setminus \{j^*\}$, given $j^*$ is of type $X$.  For $X=H, W$, let $\hat{Z}_X^S$ be the total number of individuals, including the ancestor, in the branching process $\hat{\mathcal{B}}_S$
that is defined analogously to $\hat{\mathcal{B}}_C$ in the obvious fashion.  For $X=R, H, W$ and $0 \le s_1, s_2, s_3 \le 1$, let
\begin{equation}
\label{equ:pgf_gS}
g_X^S(s_1, s_2, s_3)=\E\left[s_1^{Z_{XR}^S}s_2^{Z_{XH}^S}s_3^{Z_{XW}^S}\right].
\end{equation}
Then, cf.~\eqref{equ:progenypgfhat}, for each $s \in [0,1]$, the PGFs $f_{\hat{Z}_H^S}(s)$ and $f_{\hat{Z}_W^S}(s)$ of ${\hat{Z}_H^S}$ and ${\hat{Z}_W^S}$, respectively, are determined by
\begin{equation}
\label{equ:progenypgfhatS}
f_{\hat{Z}_H^S}(s)=s g_H^S(s, f_{\hat{Z}_H^S}(s), f_{\hat{Z}_W^S}(s)) \qquad \text{and} \qquad f_{\hat{Z}_W^S}(s)=s g_W^S(s, f_{\hat{Z}_H^S}(s), f_{\hat{Z}_W^S}(s)).
\end{equation}
Let $S$ be the total progeny of $\mcB_S$.  Then, cf.~\eqref{equ:progenypgf},
\begin{equation}
\label{equ:progenypgfS}
f_S(s)=(1-\theta)sg_R^S(s, f_{\hat{Z}_H^S}(s), f_{\hat{Z}_W^S}(s))+\theta g_H^S(s, f_{\hat{Z}_H^S}(s), f_{\hat{Z}_W^S}(s))f_{\hat{Z}_W^S}(s).
\end{equation}

Recall from~\eqref{equ:Rstar} that $R_*=\beta_G \mu_S$, since~\eqref{equ:ECnequalsESn} implies $\mu_C=\mu_S$.  Conditioning on the first generation of $\hat{\mathcal{B}}_S$ or suitable differentiation of~\eqref{equ:progenypgfhatS} yields
\[
\begin{pmatrix}
\mu_{\hat{Z}_H^S}\\
\mu_{\hat{Z}_W^S}
\end{pmatrix}
=
\begin{pmatrix}
1+\mu_{Z_{HR}^S}\\
1+\mu_{Z_{WR}^S}
\end{pmatrix}
+
\mathbf{M}^S
\begin{pmatrix}
\mu_{\hat{Z}_H^S}\\
\mu_{\hat{Z}_W^S}
\end{pmatrix},
\]
where
\begin{equation} \label{eq:M_S}
\mathbf{M}^S = \begin{pmatrix} \mu_{Z_{HH}^S} & \mu_{Z_{HW}^S}  \\
 \mu_{Z_{WH}^S}& \mu_{Z_{WW}^S} \end{pmatrix}.
\end{equation}
Let $\zeta^S$ be the maximum eigenvalue (Perron-Frobenius root) of $\mathbf{M}^S$.  In Section~\ref{sec-pig0-genI} we show that $\zeta^S$ serves as a threshold parameter for the epidemic with $\pi_G=0$, i.e.~for the epidemic with only local infection, so we also denote $\zeta^S$ by $R_L$.  Thus,
\begin{equation}
\label{equ:zetaS}
R_L (= \zeta^S) =\frac{\mu_{Z_{HH}^S}+\mu_{Z_{WW}^S}+\sqrt{(\mu_{Z_{HH}^S}-\mu_{Z_{WW}^S})^2+4 \mu_{Z_{HW}^S} \mu_{Z_{WH}^S}}}{2}.
\end{equation}
Observe that $\mu_{\hat{Z}_H^S}=\mu_{\hat{Z}_W^S}=\infty$ if $R_L \ge 1$, otherwise
\[
\begin{pmatrix}
\mu_{\hat{Z}_H^S}\\
\mu_{\hat{Z}_W^S}
\end{pmatrix}
=
\left(I-\mathbf{M}^S\right)^{-1}
\begin{pmatrix}
1+\mu_{Z_{HR}^S}\\
1+\mu_{Z_{WR}^S}
\end{pmatrix}.
\]
Differentiation of~\eqref{equ:progenypgfS} then yields
\begin{equation}
\label{equ:RstarS}
R_*=
\begin{cases}
\beta_G \left[ 1-2\theta +(1-\theta)(\mu_{Z_{RR}^S}+\mu_{Z_{RH}^S} \mu_{\hat{Z}_H^S}+\mu_{Z_{RW}^S} \mu_{\hat{Z}_W^S})+\theta(\mu_{\hat{Z}_H^S}+\mu_{\hat{Z}_W^S}) \right] & \text{ if } R_L < 1 ,\\
\infty & \text { if } R_L \ge 1.
\end{cases}
\end{equation}

\subsection{Numerical computation}
\label{subsec-numcomp}
We outline computation of $R_*$ and $z$ (and hence also $\rho$ when $I \equiv 1$).  Computation of $\rho$ when $I$ is non-constant is indicated in Section~\ref{subsec-outbreakprob} above but its implementation is computationally more intensive and intricate, owing to the appreciably larger dimension of the approximation branching process $\mathcal{B}_C$.

For $X=R, H, W$, let $\bZ_X^S=(Z_{XR}^S, Z_{XH}^S,Z_{XW}^S)$.  Note that $\bZ_X^S$ has finite support, so its joint PGF  $g_X^S(s_1, s_2, s_3)$ and mean vector are easily computed from its joint PMF (probability mass function).  Thus, given the joint PMFs of $\bZ_X^S$ $(X=R ,H, W)$, the PGF $f_S(s)$ is easily computed using~\eqref{equ:progenypgfhatS} and~\eqref{equ:progenypgfS}, enabling the limiting final size $z$ to be computed using~\eqref{equ:z}, and $R_*$ is easily computed using~\eqref{equ:RstarS}.  We outline how the joint PMF of $\bZ_R^S$ may be computed; the joint PMFs of $\bZ_H^S$ and $\bZ_W^S$ may be computed similarly.

Recall the binomially distributed independent random variables $M_1, M_2, \dots, M_d$ defined in Section~\ref{subsubsec-clump}.  Consider the susceptibility set, $\mathcal{S}_{C,j^*}$ say, of a single individual $j^*$ within a typical complex. Given $j^*$ is a remainer, the number of individuals in the complex in groups $2j-1$ and $2j$ $(j=1,2,\dots, d)$ are $n-M_j$ and $M_j$, respectively, and the number of individuals in group $2d+1$ is $\sum_{i=1}^d M_i$.  Let $\bS^C_R=(S^C_{R,1},S^C_{R,2}, \dots, S^C_{R,2d+1})$, where $S^C_{R,i}$ is the number of group-$i$ individuals in the complex susceptibility set $\mathcal{S}_{C,j^*}$, excluding $j^*$ if $i=1$.
Let $\bM=(M_1, M_2, \dots, M_d)$ and $\bm$ denote a realisation of $\bM$. For each possible $\bm$, the joint PMF of $\bS^C_R|\bM=\bm$ may be computed using Ball~\cite{Ball19}, Lemma 4.1. The joint PMF of $\bS^C_R$ may then be computed using the law of total probability and the joint PMF of $\bZ_R^S$ follows easily.

Let $\bell=(\ell_1, \ell_2, \dots, \ell_{2d+1})$ be a possible realisation of $\bS^C_R$.  Using Ball~\cite{Ball19}, Lemma 4.1, $\P(\bS^C_R=\bell|\bM=\bm)$ takes the form
\[
\P(\bS^C_R=\bell|\bM=\bm)=\alpha_R(\bm,\bell) \beta_R(\bell),
\]
where $\beta_R(\bell)$ is defined in terms of multivariate Gontcharoff polynomials (Lef\`evre and Picard~\cite{LP90}).   The $\alpha_R(\bm,\bell)$ are readily computed but the $\beta_R(\bell)$ need to be obtained by solving a triangular system of linear equations, which is numerically prohibitive unless $d$ and $n$ are suitably small.  However, note that the $\beta_R(\bell)$ are independent of $\bm$ and hence need only computing once.

\subsection{Proof of Theorem~\ref{prop-pmajor} when $\pi_G>0$}
\label{subsec-majorprobproof}
We prove here a version of Theorem~\ref{prop-pmajor} in which $\log n$ is replaced by $\kappa_n=\log \log n$.  The result for $\log n $ follows from the proof of Theorem~\ref{prop-finsize} in Section~\ref{subsec-finsizeproof}.  We use the construction of the final size of $\En$ via local infectious clumps outlined in Section~\ref{sec-genframework}.  Let $\Zn$ denote the final size of $\En$ and $\Znch$ be the final size of a modified epidemic, $\Ench$, in which the size of local infectious clumps is capped at $\kappa_n+1$.  The latter can be achieved by growing the local infectious clumps in real time and stopping the growth of a clump if and when its size reaches $\kappa_n+1$.  Thus, in $\Ench$, the local infectious clump of individual $i$ is $\mcCnch_i$, where $\mcCnch_i \subseteq \Cnmc_i$ and $\mcCnch_i$ has size
\[
\Cnch_i=
\begin{cases}
\Cn_i & \text{ if } \Cn_i \le \kappa_n+1,\\
\floor{\kappa_n}& \text { if }  \Cn_i > \kappa_n+1.
\end{cases}
\]

Let $\mcBn$ be the branching process of local infected clumps defined analogously to the branching process $\mcB$ in Section~\ref{sec-genframework}, except $(C, A, R)$ is distributed as $(\Cn, \An, \Rn)$.  Let $\Zn_I$ be the total number of infected individuals in $\mcBn$, i.e.~the sum of the sizes of the local infectious clumps of all individuals in $\mcBn$.  Let $\mcBnch$ be a modified version of $\mcBn$
in which the size of local infectious clumps is capped at $\kappa_n$ and let $\Znch_I$ be the total number of infected individuals in $\mcBnch$.  Observe that
\begin{equation}
\label{equ:ZnZnch}
\P(\Zn \le \kappa_n)=\P(\Znch \le \kappa_n) \qquad\text{and}\qquad \P(\Zn_I \le \kappa_n)=\P(\Znch_I \le \kappa_n).
\end{equation}

We construct coupled realisations of $\Ench$ and $\mcBnch$ as follows.  Let $\Gnt$ be the random graph with vertex set $\mcN$ in which for distinct $i, j\in \mcN$ there is an edge between $i$ and $j$ if and only if $i$ and $j$ belong to the same complex.  For $(i,j) \in \mcN^2$, let $\rho(i,j)$ be the length of the shortest path between $i$ and $j$ in $\Gnt$, with $\rho(i,i)=0$.  Let $\Jn_1, \Jn_2, \dots$ be i.i.d., each uniformly distributed on $\mcN$.  The initial infective in $\Ench$ is individual $\Jn_1$ and for $k=1,2, \dots$, the $k^{\rm th}$ global contact in $\Ench$ is with individual $\Jn_k$.  Let
\[
\Fn=\bigcap_{i=1}^{\floor{\kappa_n}-1}\bigcap_{j=i+1}^{\floor{\kappa_n}} \{\rho(\Jn_i,\Jn_j)>2 (\kappa_n+1)\}.
\]
If $\Fn$ occurs, then $\mcCnch_{\Jn_k}$ $(k=1,2, \dots, \floor{\kappa_n})$ are independent, so $(\Cnch_{\Jn_k},\Rnch_{\Jn_k})$ $(k=1,2, \dots, \floor{\kappa_n})$ can be used to give the clump sizes and offsprings of the first $\floor{\kappa_n}$ individuals in $\mcBnch$.  (Here, $\Rnch_i$ is the number of global contacts that emanate from $\mcCnch_i$.)  If $\Fn^C$ occurs, then we construct a realisation of $\mcBnch$ independently of $\Ench$.  Now,
\[
\P(\Znch \le \kappa_n)=\P(\Znch \le \kappa_n, \Fn)+\P(\Znch \le \kappa_n, \Fn^C)
\]
and
\[
\P(\Znch_I \le \kappa_n)=\P(\Znch_I \le \kappa_n, \Fn)+\P(\Znch_I \le \kappa_n, \Fn^C).
\]
By construction, if $\Fn$ occurs, then $\Znch \le \kappa_n$ if and only if $\Znch_I \le \kappa_n$, so
\[
\left|\P(\Znch \le \kappa_n)-\P(\Znch_I \le \kappa_n)\right| \le \P(\Fn^C).
\]

Note that, since each individual belongs to 1 or 2 complexes and each complex contains at most $2w$ individuals, for $K=1,2,\dots$, there are at most $H(K)=2\sum_{i=1}^K (2w)^i \le 4(2w)^K$ individuals within $\rho$-distance $K$ of a given individual, so
\[
\P(\Fn^C) \le \frac{\kappa_n(\kappa_n-1)4(2w)^{2(\kappa_n+1)}}{n} \to 0 \quad \text{as }n \to \infty.
\]
Thus, recalling~\eqref{equ:ZnZnch},
\begin{equation}
\label{equ:PZnZnIdiff}
\lim_{n \to \infty} \left|\P(\Zn \le \kappa_n)-\P(\Zn_I \le \kappa_n)\right|=0.
\end{equation}

Recall that $\xi$ is the extinction probability of $\mcB$ and let $\xi_n$ be the extinction probability of $\mcBn$.  Now $\Rn \convD R$ as $n \to \infty$, so $\lim_{n \to \infty}\xi_n=\xi$ (cf.~Britton et al.~\cite{BJM07}, Lemma 4.1).  Hence,
\begin{equation}
\label{equ:limsup}
\limsup_{n \to \infty} \P(\Zn_I \le \kappa_n) \le \limsup_{n \to \infty} \P(\Zn_I<\infty)=\limsup_{n \to \infty}\xi_n=\xi.
\end{equation}
Also, for fixed $k=1,2,\dots$,
\begin{equation}
\label{equ:liminf}
\liminf_{n \to \infty} \P(\Zn_I \le \kappa_n) \ge \liminf_{n \to \infty} \P(\Zn_I \le k)=\P(Z_I \le k),
\end{equation}
where $Z_I$ is the total number of individuals infected in $\mcB$.
Now, $\P(Z_I \le k) \to \xi$ as $k \to \infty$, so letting $k \uparrow \infty$ in~\eqref{equ:liminf} yields
\[
\liminf_{n \to \infty}\P(\Zn_I \le \kappa_n) \ge \xi,
\]
which together with~\eqref{equ:limsup} implies
$\lim_{n \to \infty} \P(\Zn_I \le \kappa_n)= \xi$.  Thus, using~\eqref{equ:PZnZnIdiff},
\[
\lim_{n \to \infty} \P(\Zn \le \kappa_n)=\xi,
\]
whence $\lim_{n \to \infty} \P(\Zn < \kappa_n)=\rho=1-\xi$,
which completes the proof of Theorem~\ref{prop-pmajor}, with $\log n$ replaced by $\kappa_n$.

Note that a similar argument to the above shows that
\[
\lim_{n \to \infty}\P(\Zn \le k)=\P(Z_I \le k)\qquad(k=1,2,\dots),
\]
so $\Zn \convD Z_I$ as $n \to \infty$.  The random variable $Z_I$ is defective, with a mass at infinity, if $R_*>1$.

\subsection{Proof of Theorem~\ref{prop-finsize} when $\pi_G>0$}
\label{subsec-finsizeproof}

To prove Theorem~\ref{prop-finsize}, we use the adaptation of the embedding technique of Scalia-Tomba~\cite{ScaliaTomba85, ScaliaTomba90} to the general two-levels of mixing setting; see Ball and Neal~\cite{BN02}, Section 6.2 and Ball and Sirl~\cite{BS19}, Section 3.4.2.

Let $I_1, I_2, \dots, I_n$ and the random directed graph $\Gn$ be as in Section~\ref{sec-genframework} and let $L_1, L_2, \dots, L_n$ be an independent sequence of i.i.d. exponential random variables, each having rate $\beta_G$.  Given, $I_1, I_2, \dots, I_n$, $\Gn$ and $L_1, L_2, \dots, L_n$, we construct for each $t \ge 0$ a realisation of the final outcome of an epidemic $\Ent(t)$ as follows.  Note that in this construction $t$ refers to global infectious pressure and not time.  We determine first who is infected globally; specifically, for $i \in \mathcal{N}$, individual $i$ is infected globally if and only if $L_i \le t$.  We then use the random directed graph $\Gn$ to determine who is subsequently infected locally.  Specifically, for $i \in \mathcal{N}$, individual $i$ is infected locally if and only if at least one member of $i$'s local susceptibility set $\Snmc_i$ is infected globally.
For $i \in \mcN$, let
\begin{align*}
\chin_i(t)&=1_{\{i \text{ is infected in } \Ent(t)\}}\\
&=
\begin{cases}
1 & \text{ if } \min_{j \in \Snmc_i} L_j \le t,\\
0 & \text { otherwise } .
\end{cases}
\end{align*}
Let
\[
\rbn(t)= \sum_{i=1}^n \chin_i(t)\qquad\text{and}\qquad \abn(t)= \sum_{i=1}^n I_i \chin_i(t)
\]
be respectively the total size and severity of the epidemic $\Ent(t)$.

To connect with the epidemic $\En$, consider an epidemic, $\hat{\mathcal{E}}^{(n)}$ say, that is initiated by exposing the population to $\Tn_0$ units of global infectious pressure, so each individual is exposed to $\bartn_0=n^{-1}\Tn_0$ units of global infectious pressure.  These $\Tn_0$ units of global infectious pressure may infect some individuals globally and the consequent local spread will lead to $\abn(\bartn_0)$ further units of global infectious pressure. Thus, the population is now exposed to $\Tn_0+\abn(\bartn_0)$ units of global infectious pressure and the process can be continued in the obvious fashion. For $k=0,1,\dots$, let
\[
\bartn_{k+1}=\bartn_0+\abn(\bartn_k).
\]
Then $k_*^{(n)}=\min\{k:\bartn_{k+1}=\bartn_k\}$ is well-defined as the population is finite.  Let $\Tn_{\infty}=\Tn_{k_*^{(n)}}$ and $\bartn_{\infty}=n^{-1}\Tn_{\infty}$. Note that
\begin{equation}
\label{equ:tbarinf}
\bartn_{\infty}=\min\{t \ge 0: t=\bartn_0+\barzan(t)\},
\end{equation}
where $\barzan(t)=n^{-1}\abn(t)$.  The total size and severity of the epidemic $\hat{\mathcal{E}}^{(n)}$ are given by $\rbn(\bartn_{\infty})$ and $\abn(\bartn_{\infty})$, respectively.

The epidemics $\En$ and $\hat{\mathcal{E}}^{(n)}$ have different laws, since $\En$ is initiated by a single individual chosen uniformly at random from the population being externally infected and the number of initial infectives in $\hat{\mathcal{E}}^{(n)}$, $\hat{Z}^{(n)}_{\rm init}$ say, is distributed as $\rbn(\bartn_0)$. However, owing to the lack-of-memory property of the exponential distribution, $\En$ has the same law as $\hat{\mathcal{E}}^{(n)}$ conditioned upon $\hat{Z}^{(n)}_{\rm init}=1$.  Suppose that $\P(\Tn_0=1)=1$.  Then, $\hat{Z}^{(n)}_{\rm init} \sim \text{Bin}(n, 1-\exp(-\beta_G/n))$, so
\[
\P(\hat{Z}^{(n)}_{\rm init}=1)=n \exp(-(n-1)\beta_G/n)(1-\exp(-\beta_G/n)) \to \beta_G \exp(-\beta_G) \quad \text{as } n \to \infty.
\]
Thus, $\liminf_{n \to \infty}\P(\hat{Z}^{(n)}_{\rm init}=1)>0$ and (cf.~Jansen~\cite{Jansen09}) it follows that convergence in probability results for $\hat{\mathcal{E}}^{(n)}$ hold also for $\hat{\mathcal{E}}^{(n)}|\hat{Z}^{(n)}_{\rm init}=1$, and hence for $\En$.

Recall the random variable $S$ defined at~\eqref{equ:CSAconv} in Section~\ref{sec-genframework}, describing the size of a typical local susceptibility set in the limit $n \to \infty$.  For $t \ge 0$, let
$r(t)=1-f_S(\re^{-\beta_G t})$ and $\barzrn(t)=n^{-1}\rbn(t)$.

${\newtheorem{lemma-glivenko}[theorem]{Lemma}}$
\begin{lemma-glivenko} \label{lem:glivenko}
Suppose that ${\var}(I)<\infty$.  Then, for any $t_0>0$,
\begin{equation}
\label{equ:glivenko}
\sup_{t_0 \le t < \infty} |\barzrn(t)-r(t)| \convp 0 \quad\text{and}\quad \sup_{t_0 \le t < \infty} |\barzan(t)-r(t)| \convp 0 \quad \text{as } n \to \infty.
\end{equation}
\end{lemma-glivenko}

\begin{proof}
We prove the result for $\barzan(t)$.  The result for $\barzrn(t)$ is proved similarly.  Recall that we assume $\E[I]=1$.

For $t \ge 0$, let
\begin{align*}
\an(t)=\E[\barzan(t)]&=\E[I_1 \chin_1(t)]\\
&=\E[I_1(1-\exp(-\beta_G t \Sn_1))]\\
&=1-\E[\exp(-\beta_G t \Sn_1)],
\end{align*}
since $I_1$ and $\Sn_1$ are independent.  We show below that, for each $t \in (0, \infty)$, $\barzan(t)-\an(t) \convp 0$ as $n \to \infty$.  Now $\an(t) \to r(t)$ as $n \to \infty$, since $\Sn_1 \convD S$ as $n \to \infty$ and $0 \le \exp(-\beta_G t \Sn_1) \le 1$ for all $n$.  Thus, $\barzan(t) \convp r(t)$ as $n \to \infty$.  This convergence holds also for $t=\infty$, since $\barzan(\infty)=n^{-1}\sum_{i=1}^n I_i \convp 1=r(\infty)$ as $n \to \infty$ by the weak law of large numbers.  The random function $\barzan(t)$ is non-decreasing in $t$ on $[t_0, \infty]$ and $r(t)$ is continuous and increasing on $[t_0,\infty]$. The second result in~\eqref{equ:glivenko} then follows by a similar argument to the proof of Ball and Lyne~\cite{BL01}, Lemma 5.1, see Ball and Sirl~\cite{BS19}, page 187.

To complete the proof, we show that for each $t \in (0, \infty)$, $\barzan(t)-\an(t) \convp 0$ as $n \to \infty$.
By Chebysvev's inequality
\[
\P(|\barzan(t)-\an(t)|> \epsilon) \le \frac{\var(\barzan(t))}{\epsilon^2}=\frac{\var(\abn(t))}{n^2 \epsilon^2},
\]
so a sufficient condition for $\barzan(t)-\an(t) \convp 0$ as $n \to \infty$ is
\begin{equation}
\label{equ:varantcond}
\lim_{n \to \infty} n^{-2}\var(\abn(t))=0.
\end{equation}

For $t \ge 0$, let $\chint_i(t)=1-\chin_i(t)$ and $\abnt(t)=\sum_{i=1}^n I_i \chint_i(t)$.  Let $\sigma_I^2=\var(I)$.  Noting that $\abn(t)+\abnt(t)=\sum_{i=1}^n I_i$, we have that
\[
n\sigma_I^2=\var(\abn(t))+\var(\abnt(t))+2\cov(\abn(t),\abnt(t)).
\]
Further,
\begin{align*}
\cov(\abn(t),\abnt(t))&=\cov\left(\sum_{i=1}^n I_i(1-\chint_i(t)), \sum_{i=1}^n I_i \chint_i(t) \right)\\
&= \cov\left(\sum_{i=1}^n I_i, \sum_{i=1}^n I_i \chint_i(t)\right)-\var(\abnt(t)),
\end{align*}
so~\eqref{equ:varantcond} holds if
\begin{equation}
\label{equ:varantcond0}
\lim_{n \to \infty} n^{-2} \var(\abnt (t))=0\qquad \text{and}\qquad \lim_{n \to \infty} n^{-2} \cov\left(\sum_{i=1}^n I_i, \sum_{i=1}^n I_i \chint_i(t)\right)=0.
\end{equation}

For $t \ge 0$, let $\abnh(t)=\sum_{i=1}^n I_i \chint_i(t)1_{\{\Sn_i \le \kappa_n\}}$, where $\kappa_n=\log \log n$ as in Section~\ref{subsec-majorprobproof}.  Then, using the independence of $I_i$ and $\Sn_i$ and recalling that $\E[I]=1$,
\begin{align}
\label{equ:meandiff}
\left|\E[\abnt(t)]-\E[\abnh(t)]\right|&=\E\left[\sum_{i=1}^n I_i \chint_i(t)1_{\{\Sn_i > \kappa_n\}}\right]\nonumber\\
&=\sum_{i=1}^n \E \left[I_i \exp(-\beta_G t \Sn_i)1_{\{\Sn_i > \kappa_n\}}\right]\nonumber\\
&\le ng_n(t),
\end{align}
where $ g_n(t) = \exp(-\beta_G t \kappa_n)$.
Further,
\begin{align*}
\E\left[\abnt(t)^2\right]&=\E\left[\sum_{i=1}^n I_i \chint_i(t) \sum_{j=1}^n I_j \chint_i(t)\right]\\
&=\sum_{i=1}^n \sum_{j=1}^n \E\left[I_i I_j \exp(-\beta_G t |\Snmc_i \cup \Snmc_j|)\right].
\end{align*}
Thus,
\begin{align}
\label{equ:meandiff1}
\left| \E[\abnt(t)^2]-\E[\abnh(t)^2]\right|&=\sum_{i=1}^n \sum_{j=1}^n\E\left[I_iI_j\exp(-\beta_G t |\Snmc_i \cup \Snmc_j|)\left(1-1_{\{\Sn_i \le \kappa_n\}}1_{\{\Sn_j \le \kappa_n\}}\right)\right]\nonumber\\
&\le \sum_{i=1}^n \sum_{j=1}^n\E\left[I_iI_j\exp(-\beta_G t |\Snmc_i \cup \Snmc_j|)\left(1_{\{\Sn_i > \kappa_n\}}+1_{\{\Sn_j > \kappa_n\}}\right)\right]\nonumber\\
&\le 2g_n(t)\sum_{i=1}^n \sum_{j=1}^n \E[I_i I_j]\nonumber\\
&=2g_n(t)[n\sigma_I^2+n^2].
\end{align}
It follows from~\eqref{equ:meandiff} and~\eqref{equ:meandiff1} that, for $t>0$,
\begin{equation}
\label{equ:varantcond1}
\lim_{n \to \infty} n^{-2} \var(\abnt (t))=0 \iff \lim_{n \to \infty} n^{-2} \var(\abnh (t))=0,
\end{equation}
since $g_n(t) \to 0$ as $n \to \infty$.

Note that, by exchangeability,
\begin{align}
\label{equ:varabhat}
n^{-2} \var(\abnh (t))&=n^{-2} \sum_{i=1}^n \sum_{j=1}^n \cov\left(I_i \chint_i(t)1_{\{\Sn_i \le \kappa_n\}}, I_j \chint_j(t)1_{\{\Sn_j \le \kappa_n\}} \right)\nonumber\\
&=n^{-1}\sum_{j=1}^n \cov\left(I_1 \chint_1(t)1_{\{\Sn_1 \le \kappa_n\}}, I_j \chint_j(t)1_{\{\Sn_j \le \kappa_n\}} \right)\nonumber\\
&=\cov\left(I_1 \chint_1(t)1_{\{\Sn_1 \le \kappa_n\}}, I_J \chint_J(t)1_{\{\Sn_J \le \kappa_n\}} \right),
\end{align}
where $J$ is chosen uniformly at random from $\mcN$.

Recall from Section~\ref{subsec-majorprobproof} the graph $\Gnt$ and the associated distance function $\rho(i,j)$.  Let $\Xn_i(t)=I_i \chint_i(t)1_{\{\Sn_i \le \kappa_n\}}$ $(i \in \mcN)$ and $\Un=1_{\{\rho(1,J)\le 3 \kappa_n\}}$.  Then,
\begin{align*}
\cov\left(\Xn_1(t),\Xn_J(t)\right)&=\E\left[\cov\left( \Xn_1(t),\Xn_J(t) \mid \Un\right)\right]\\
&\qquad\qquad+\cov\left(\E\left[\Xn_1(t)\mid \Un \right], \E\left[\Xn_J(t)\mid \Un \right]\right).
\end{align*}
Let $p_n=\P(\Un=1)$ and note that $p_n \le \frac{H(3\kappa_n)}{n}$, where $H(K)\le 4(2w)^K$; see just before~\eqref{equ:PZnZnIdiff} in Section~\ref{subsec-majorprobproof}.  Hence, $p_n \to 0$ as $n \to \infty$. Further, $\Xn_1(t)$ and $\Xn_J(t)$ are independent if $\Un=0$, since $\{i \in \mcN: \rho(1,i) \le \kappa_n+1\}$ and  $\{j \in \mcN: \rho(j,J) \le \kappa_n+1\}$ are disjoint.  Thus,
\begin{align*}
\left|\E\left[\cov\left( \Xn_1(t),\Xn_J(t) \mid \Un\right)\right]\right| &\le p_n \left|\cov\left( \Xn_1(t),\Xn_J(t) \mid \Un=1\right)\right|\\
&\le p_n \E[I^2] \to 0\quad \text{as } n \to \infty.
\end{align*}
Also,
\begin{align*}
\left|\cov\left(\E\left[\Xn_1(t)\mid \Un \right],\right.\right.&\left.\left. \E\left[\Xn_J(t)\mid \Un \right]\right)\right|\\
&=p_n(1-p_n)\left(\E\left[\Xn_J(t)\mid \Un=1 \right]-\E\left[\Xn_J(t)\mid \Un=0 \right]\right)^2\\
&\le p_n(1-p_n)\E[I]^2 \to 0\quad \text{as } n \to \infty.
\end{align*}
Thus, $\cov\left(\Xn_1(t),\Xn_J(t)\right) \to 0$ as $n \to \infty$, so the first condition in~\eqref{equ:varantcond0} is satisfied using~\eqref{equ:varantcond1} and~\eqref{equ:varabhat}. The second condition in~\eqref{equ:varantcond0} is proved similarly, so the details are left to the reader.
\end{proof}
${\newtheorem{remark1}[theorem]{Remark}}$
\begin{remark1}
Note that $\barzan(0)=\an(0)=0$, so~\eqref{equ:glivenko} holds with $t_0=0$ if $\lim_{t \downarrow 0} r(t) = 0$, a necessary and sufficient condition for which is $R_L \le 1$.
\end{remark1}

Suppose that $\P(\Tn_0=1)=1$ and that $R_*>1$. Note that the equation $r(t)=t$ is the same equation satisfied by $z$ at~\eqref{equ:z} in Section~\ref{sec-genframework}.  Now, $r(0)=0$, $r'(0)=R_*>1$ and $r''(t)<0$ for $t \ge 0$.  Thus, $r(t)$ is concave on $[0, \infty)$.  Further, $r(\infty)=1$, so there exists a unique $z \in (0,\infty)$ satisfying $r(z)=z$ and $r'(z)<1$.  (If $r'(z)=1$ then $r(z)=\int_0^{z} r'(x)\,{\rm d}x> \int_0^{z}1 \,{\rm d}x=z$, contradicting $r(z)=z$.)  It follows from~\eqref{equ:tbarinf} and the second result in~\eqref{equ:glivenko} that
$\min(\bartn_{\infty}, |\bartn_{\infty}-z|) \convp 0$ as $n \to \infty$, and hence from the first result in~\eqref{equ:glivenko} that
\begin{equation}
\label{equ:barzlim}
\min(\barzrn(\bartn_{\infty}),|\barzrn(\bartn_{\infty})-z|) \convp 0 \quad\text{as } n \to \infty.
\end{equation}

Recall that $\barzrn(\bartn_{\infty})=n^{-1}\rbn(\bartn_{\infty})$, where $\rbn(\bartn_{\infty})$ is the total size of the epidemic $\hat{\mathcal{E}}^{(n)}$, and that convergence in probability results that hold for $\hat{\mathcal{E}}^{(n)}$, hold also for $\En$.  Let $\Znb=n^{-1}\Zn$, where $\Zn$ is the total size of the epidemic $\En$.  It follows from the above that $\Znb \convp \bar{Z}$ as $n \to \infty$, where the random variable $\bar{Z}$ has support $\{0, z\}$.  All that remains is to determine $\P(\bar{Z}=0)$.  Using Theorem~\ref{prop-pmajor}, with $\log n$ replaced by $\kappa_n$,
\begin{equation}
\label{equ:liminf1}
\liminf_{n \to \infty} \P(\Znb \le \frac{z}{2}) \ge \liminf_{n \to \infty} \P(\Zn \le \kappa_n)=\xi.
\end{equation}

To obtain an upper bound, first note using~\eqref{equ:barzlim} that, for any $\epsilon \in (0,z/2)$,
\begin{equation}
\label{equ:limsup1}
\limsup_{n \to \infty}\P(\Znb \le \frac{z}{2})= \limsup_{n \to \infty}\P(\Znb \le \epsilon).
\end{equation}
For $K=1,2,\dots$, let $\Ench_K$ be a modification of the epidemic $\En$ in which the size of all local infectious clumps is capped at $K$ and let $\Znch_K$ be the total size of $\Ench_K$.  Clearly,
\begin{equation}
\label{equ:ZnZnKbound}
\P(\Znb \le \epsilon) \le \P(\Znch_K \le n \epsilon).
\end{equation}
To derive an upper bound for $\P(\Znch_K \le n \epsilon)$, we extend a technique which dates back to Whittle~\cite{Whittle55}.  A realisation of
$\Ench_K$ can constructed in an analogous fashion to that of $\Ench$ in Section~\ref{subsec-majorprobproof}, by using a sequence $\Jn_1, \Jn_2, \dots$ of i.i.d. random variables that are uniformly distributed on $\mcN$ to determine the individuals contacted by successive global contacts.  A lower bound process may be obtained by deleting any global contact (and ignoring the consequent local infectious clump) that is with an individual whose $\rho$-distance from any previously globally contacted individual (including the initial infective) is $\le 2K$.  While the total size of $\Ench_K$ is $\le \epsilon n$, the probability that a global contact is deleted is bounded above by $\frac{\epsilon n H(2K)}{n} \le 4\epsilon (2w)^{2K}$.  Thus,
\begin{equation}
\label{equ:ZnKZnIKbound}
\P(\Znch_K \le n \epsilon) \le \P(\Znch_{I,K, \epsilon} \le n \epsilon),
\end{equation}
where $\Znch_{I,K, \epsilon}$ is the total number of individuals infected in the branching process $\mcBnch_{K, \epsilon}$, whose law is similar to that of $\mcBn$, except the size of local infectious clumps is capped at $K$ and each individual (i.e.~local infectious clump) is deleted at birth independently with probability $\min(4\epsilon (2w)^{2K},1)$.  Let $\mcB_{K, \epsilon}$ be an analogous branching process derived from $\mcB$.  Let $\xi_n(K, \epsilon)$ and $\xi(K, \epsilon)$ be the extinction probabilities of $\mcBnch_{K, \epsilon}$ and $\mcB_{K, \epsilon}$, respectively, and note that $\lim_{n \to \infty} \xi_n(K, \epsilon)=\xi(K, \epsilon)$.  Using~\eqref{equ:limsup1}, \eqref{equ:ZnZnKbound} and~\eqref{equ:ZnKZnIKbound},
\begin{equation}
\label{equ:limsup2}
\limsup_{n \to \infty}\P(\Znb \le \frac{z}{2})\le \limsup_{n \to \infty} \P(\Znch_{I,K, \epsilon} \le n \epsilon)\\
\le \limsup_{n \to \infty} \xi_n(K, \epsilon) =\xi(K, \epsilon).
\end{equation}
Now $\xi(K, \epsilon) \downarrow \xi(K)$ as $\epsilon \downarrow 0$, where $\xi(K)$ is the extinction probability for the branching process derived from $\mcB$ in which clumps are capped at size $K$, and $\xi(K) \downarrow \xi$ as $ K \uparrow \infty$.  Thus,~\eqref{equ:liminf1} and~\eqref{equ:limsup2} yield
\[
\lim_{n \to \infty}\P(\Znb \le \frac{z}{2})=\xi,
\]
whence $\P(\bar{Z}=0)=\xi=1-\P(\bar{Z}=z)$.

Recall from the end of Section~\ref{subsec-majorprobproof} that $\Zn \convD Z_I$ as $n \to \infty$. It follows easily from the above that, if $(\upsilon_n)$ is any sequence of nonnegative integers satisfying
$\upsilon_n \to \infty$ and $n^{-1}\upsilon_n \to 0$ as $n \to \infty$, then
\[
\lim_{n \to \infty}\P(\Zn \ge \upsilon_n)=1-\xi
\]
and, if $R_*>1$,
\[
\Znb|\Zn \ge \upsilon_n \convp z \quad \text{as } n \to \infty.
\]

\section{Derivations and proofs when $\pi_G=0$} \label{sec-pig0}

\subsection{Introduction and overview} \label{sec-pig0-intro}

In this section we consider the special case where $\pi_G =0$, so that all infection is local through the households and workplaces.  
Throughout this section we assume that $\theta >0$, since otherwise the complexes and workplaces coincide and in the absence of global infection the epidemic is restricted to its original workplace.
We prove Theorems \ref{prop-pmajor} and \ref{prop-finsize} concerning the probability ($\rho$) and final size ($z$) of a major outbreak hold in this setting.  Note that $\rho=0$ unless the local infectious clumps and local susceptibility sets are supercritical, i.e.~unless $R_L>1$, in which case
$z$ and $\rho$ are given by \eqref{eq:pig0:extinctS} and \eqref{eq:pig0:extinct}, respectively.

Let $i^\ast$ denote the initial infective. We couple the construction of the local infections clump for individual $i^\ast$, $\Cnmc_{i^\ast}$, to a {\it forward} branching process $\mathcal{B}_C$. Then we show that
\begin{align*} 
\P (C^{(n)}> \log n) \rightarrow \rho \qquad \mbox{as } n \rightarrow \infty,
\end{align*}
where $C^{(n)}= |\mathcal{C}^{(n)}_{i^\ast}|$ and $\rho$ is the non-extinction probability of $\mathcal{B}_C$ given by \eqref{eq:pig0:extinct}. Furthermore, we show that there exists $\delta^\prime >0$ such that
\begin{align} \label{eq:pig0:rho2}
\P ( C^{(n)}> \delta^\prime n |C^{(n)}> \log n) \rightarrow 1 \qquad \mbox{as } n \rightarrow \infty.
\end{align}
In the absence of global infection, $Z^{(n)} =C^{(n)}$.

Similarly, we couple  the construction of the local infections susceptibility set, $\mathcal{S}^{(n)}$, to a {\it backward} branching process $\mathcal{B}_S$. We show that
\begin{align} \label{eq:pig0:z}
\P (S^{(n)} > \log n) \rightarrow z \qquad \mbox{as } n \rightarrow \infty,
\end{align}
where $z$ is the non-extinction probability of $\mathcal{B}_S$ given by \eqref{eq:pig0:extinctS}.
Consider a typical initially susceptible individual $j^\ast$, with local susceptibility set $\Snmc_{j^\ast}$ having size $\Sn_{j^\ast} \sim \Sn$.  If $\Sn_{j^\ast} < \log n$, then with high probability (i.e.~with probability tending to one as $n \to \infty$), $\Snmc_{j^\ast}$ does not contain $i^\ast$, so $j^\ast$ is not infected by the epidemic.  Suppose instead that $\Sn_{j^\ast} \ge \log n$.  Then, in the event of a major outbreak, it follows from ~\eqref{eq:pig0:rho2} that $C^{(n)}> \delta^\prime n$ with high probability, and consequently $\Snmc_{j^\ast} \cap \Cnmc$ is non-empty with high probability and $j^\ast$ is infected by the epidemic.  The events $\Sn_{j_1^\ast} < \log n$ and $\Sn_{j_2^\ast} < \log n$, where $j_1^\ast$ and $j_2^\ast$ are chosen uniformly at random from the initial susceptibles, are asymptotically independent and consequently $n^{-1} \Zn | Z^{(n)} > \log n \convp z$ as $n \to \infty$. This is the intuition underlying the proof.  Rather than considering whether $\Sn_{j^\ast} < \log n$, it is convenient to consider whether $\Snmc_{j^\ast}$ goes extinct before generation $\ell_n$, where $\ell_n = \lceil a \log \log n \rceil$ for suitable choice of $a >0$, cf.~Ball et al.~\cite{BST14}.

The remainder of this section is organised as follows. In Section~\ref{sec-pig0-clumpsus}, we derive some asymptotic properties of local infectious clumps and susceptibility sets, focusing mainly on the latter.  In Section~\ref{sec-pig0-zmajor}, we consider the final outcome of the epidemic.  We show first in Lemma~\ref{lemma:pig0:hatZ} that the fraction of the population whose local susceptibility set reaches generation $\ell_n$ converges in probability to $z$.  Then, in Lemma~\ref{lem:pig0:susset}, we consider the special case when $I \equiv 1$ and show that, conditional upon the occurrence of a major outbreak, the difference between that fraction and $n^{-1}\Zn$ converges in probability to 0 as $n \to \infty$, leading to
\begin{align*} 
\left. \frac{Z^{(n)}}{n} \right| Z^{(n)} > \log n \convp z \qquad \mbox{as } n \rightarrow \infty,
\end{align*}
i.e.~Theorem \ref{prop-finsize} in the case $\pi_G =0$ and $I \equiv 1$.  Finally, in Section~\ref{sec-pig0-genI}, we extend the proof to the case of general $I$.

\subsection{Asymptotic properties of local infectious clumps and susceptibility sets} \label{sec-pig0-clumpsus}
For some $a >0$, let $\ell_n = \lceil a \log \log n \rceil$, where further requirements on $a$ will be defined as required. We start by coupling the construction of  $\mathcal{C}^{(n)}$ to $\mathcal{B}_C$ and  $\mathcal{S}^{(n)}$ to $\mathcal{B}_S$ over the first $\ell_n$ generations, where the initial individual is generation 0.
Since the maximum size of a complex is $2 w$ (all individuals in the workplace are movers), we have that the number of individuals whose shortest path from a randomly selected individual has length at most $\ell_n$, is bounded above by
\begin{align} \label{eq:pig0:dist}
\sum_{i=0}^{\ell_n} 2 (2w)^i = 2 \frac{[2w]^{\ell_n+1} -1}{2w -1},
\end{align}
where the factor 2 comes from mover individuals belonging to two complexes. For any $\epsilon >0$, the right hand side of \eqref{eq:pig0:dist} is less than $n^\epsilon$ for all sufficiently large $n$. The total number of movers in the population is ${\rm Bin} (n, \theta)$, so the proportion of the population who are movers is approximately $\theta$. Given any single complex contains at most $2w$ movers 
it is straightforward using the birthday problem, see, for example, Ball and Donnelly \cite{BD95}, p.4, to show that with probability tending to 1, distinct complexes are used in the first $\ell_n$ generations of $\mathcal{S}^{(n)}$ or $\mathcal{C}^{(n)}$. In other words, the growth of the population from a randomly selected individual is, with high probability, a tree-like configuration of complexes over the first $\ell_n$ generations {\it c.f.}~ the proof of Lemma \ref{lem:rfclumpconv}.
This can be taken further. Let $\mathcal{K}^{(n)} = \mathcal{K}^{(n)}_S \cup \mathcal{K}^{(n)}_C$ be a set of $K = K_S + K_C \in \mathbb{N} $ randomly selected individuals in a population of size $n$, where for the $K_S$ $(K_C)$ individuals in $\mathcal{K}^{(n)}_S$ $( \mathcal{K}^{(n)}_C)$ we construct the susceptibility set (clump) to generation $\ell_n$. Then with probability tending to 1, no complex is used more than once in the construction of any of the susceptibility sets or clumps for the first $\ell_n$ generations.

In order to obtain $z$ we employ a construction of a local susceptibility set similar to the local infectious clump given in Section \ref{subsubsec-clump}. In Section \ref{subsubsec-clump} and Lemma  \ref{lem:rfclumpconv}, the focus is on the probability mass function of $C^{(n)}$ and its limit $C$. In the absence of global infection, we simply need to consider whether or not the limiting branching process, $\mathcal{B}_S$, goes extinct.

The construction of $\mathcal{B}_S$ given in Section \ref{subsubsec-susset} can again be used with $\mathbf{Z}^S_X = (Z^S_{XH}, Z^S_{XW})$ being the number of new complex susceptibility sets that are added to a local susceptibility set (offspring in the branching process) from a type $X$ complex susceptibility set.
Let $j^*$ denote the initial individual for whom we are constructing a local susceptibility set. If $j^*$ is a mover then they belong to two complexes and its offspring distribution is $\mathbf{Z}^S_H + \mathbf{Z}^S_W$, where the two random vectors are independent. If $j^*$ is a remainer, then without loss of generality they belong to group 1 and $\mathbf{Z}_R^S = (Z^S_{RH}, Z^S_{RW})$ gives the number of movers in the complex who belong to the susceptibility set of $j^*$. Note that, beyond the identity (mover or remainer) of the initial individual, we are not interested in any remainers who belong to the complex susceptibility set as they play no role in the extinction probability.

For $X=R, H, W$ and $\mathbf{s} \in [0,1]^2$, let
\begin{align*} 
\bar{g}_X^S (\mathbf{s}) = \E \left[ s_1^{Z_{XH}^S} s_2^{Z_{XW}^S} \right].
\end{align*}
${\newtheorem{remark2}[theorem]{Remark}}$
\begin{remark2}
Note that $\bar{g}_X^S (\mathbf{s}) = g_X^S (1,s_1,s_2)$ where $g_X^S$ is defined in \eqref{equ:pgf_gS}.  Thus, $\bar{g}_X^S$ can be computed using methods described in Section~\ref{subsec-numcomp}.
\end{remark2}

Let $\bmeta^S = (\eta^S_H, \eta^S_W)$ denote the smallest solution in $[0,1]^2$ of
\begin{align}  \label{eq:pig0:pgfS2}
\bmeta^S = \left( \bar{g}^S_H (\bmeta^S), \bar{g}^S_W (\bmeta^S) \right),
\end{align}
where $\eta^S_H$ $(\eta^S_W)$ is the probability that the branching process $\mathcal{B}_S$ goes extinct from a single type $H$ ($W$) individual. Then $z$, the non-extinction of $\mathcal{B}_S$, satisfies
\begin{align} \label{eq:pig0:extinctS}
1-z = (1-\theta) \bar{g}_R^S (\bmeta^S) + \theta \eta^S_H \eta^S_W.
\end{align}

For $i=1,2,\ldots$, let $\mathbf{V}_i^S = (V_{H,i}^S, V_{W,i}^S)$ denote the number of individuals in the $i^{th}$ generation of $\mathcal{B}_S$. Similarly, let  $\mathbf{V}_i^{S^{(n)}} = (V_{H,i}^{S^{(n)}} , V_{W,i}^{S^{(n)}} )$ denote the number of individuals in the $i^{th}$ generation of $\mathcal{S}^{(n)}$. Since with probability tending to 1, $\mathcal{S}^{(n)}$ can be coupled to $\mathcal{B}_S$ to coincide over the first $\ell_n$ generations, we have that
\begin{align*}  
\P \left( \mathbf{V}_{\ell_n}^{S^{(n)}}  \neq  \mathbf{V}_{\ell_n}^S\right) \rightarrow 0 \qquad \mbox{as } n \rightarrow \infty.
\end{align*}
We use properties of branching processes to explore $\mathbf{V}_{\ell_n}^S$.

The mean offspring matrix, $\mathbf{M}^S$, given in \eqref{eq:M_S}, with Perron-Frobenius eigenvalue $\zeta^S (=R_L)$, determines the behaviour of the  branching process $\mathcal{B}_S$. Note that $\bmeta^S=(1,1)$ if $\zeta^S\le 1$ and $\bmeta^S$ is the unique solution of~\eqref{eq:pig0:pgfS2} if $\zeta^S > 1$.  Hence, $z >0$ if and only if $\zeta^S >1$. Also, for $0 < \theta <1$, all elements of $\mathbf{M}^S$ are non-zero and the branching process $\mathcal{B}_S$ is aperiodic.  Therefore by Mode \cite{Mode71}, p.19, (8.2), there exists a random variable $W^S$ with $\P(W^S=0) =1-z$ and a continuous distribution on  $(0,\infty)$, such that
\begin{align} \label{eq:mode82}
\lim_{k \rightarrow \infty} \frac{1}{[\zeta^S]^k} \mathbf{V}^S_k = W^S \mathbf{v}^S \qquad \mbox{a.s.},
\end{align}
where $\mathbf{v}^S > \mathbf{0}$ is the normalised left eigenvector of $\mathbf{M}^S$ corresponding to the Perron-Frobenius root $\zeta^S$.  For $\theta=1$, the branching process $\mathcal{B}_S$ has period 2 and the arguments are easily modified by letting $\ell_n = \lceil a \log \log n \rceil+1$ if $\lceil a \log \log n \rceil$ is odd and using single-type branching process theory.

${\newtheorem{lemma-pig0-susset1}[theorem]{Lemma}}$
\begin{lemma-pig0-susset1} \label{lem:pig0:susset1}
Suppose that $\zeta^S >1$. Let $c^S = \log \zeta^S >0$ and $a= 5 /(2 c^S)$, and recall that $\ell_n = \lceil a \log \log n \rceil$. Then
 \begin{align} \label{eq:pig0:VSn2}
\P \left( [\log n]^2  < \mathbf{V}^{S^{(n)}}_{\ell_n} < [\log n]^3 \right) \rightarrow z \qquad \mbox{as } n \rightarrow \infty,
\end{align}
where $z$ satisfies \eqref{eq:pig0:extinctS}.
\end{lemma-pig0-susset1}
\begin{proof}
Note that
\[ [\zeta^S]^{\ell_n} = \exp (c^S \lceil a \log \log n \rceil) \approx \exp \left( \log [\log n]^{c^S a} \right) = [\log n]^{c^S a} . \]
Therefore
\[ \mathbf{V}^S_{\ell_n} \approx [\log n]^{c^S a} W^S \mathbf{v}^S,  \]
where $W^S=0$ if only and if the branching process goes extinct. Therefore, for any $0 < d_0 < c^S a < d_1$, conditional on non-extinction of the branching process, $\mathcal{B}_S$, we have that for all sufficiently large $n$, the number of individuals in generation $\ell_n$ lies between $[\log n]^{d_0} $ and $[\log n]^{d_1} $. Since $a c^S = 2.5$, \eqref{eq:pig0:VSn2}, follows immediately.
\end{proof}

It is straightforward to adapt the arguments used in the proof of Lemma \ref{lem:pig0:susset1} to prove \eqref{eq:pig0:z}. Let $0 < p_0 \leq 1$ then using \eqref{eq:mode82}, we have that the $\P( \mathbf{V}^S_{\lceil p_0\ell_n \rceil} \neq \mathbf{0}) \rightarrow z$ as $n \rightarrow \infty$. For sufficiently small $p_0$, the number of individuals in the first $p_0 \ell_n$ generations is less than $[\log n]$ with probability tending to 1 as $n \rightarrow \infty$. However, if the susceptibility set reaches generation $\lceil p_0\ell_n \rceil$ without going extinct, it will, with probability
tending to 1, reach generation $\ell_n$ without going extinct, as the number of individuals in generation $\lceil p_0\ell_n \rceil$ tends to infinity as $n \rightarrow \infty$.

\subsection{Size of major outbreaks} \label{sec-pig0-zmajor}

For $n=1,2,\ldots$ and $j=1,2,\ldots, n$, let $\check{S}_j^{(n)}$ denote the number of individuals in generation $\ell_n$ of the susceptibility set of individual $j$. In particular, we are interested in whether or not $\check{S}_j^{(n)} =0$, i.e.~whether or not the susceptibility set goes extinct by generation $\ell_n$. Let
\begin{align*}  
\check{Z}^{(n)} = \sum_{j=1}^n 1_{\{ \check{S}_j^{(n)} \neq 0\}}
\end{align*}
denote the number of susceptibility sets which do not go extinct by generation $\ell_n$.
${\newtheorem{lemma-pig0-hatZ}[theorem]{Lemma}}$
\begin{lemma-pig0-hatZ} \label{lemma:pig0:hatZ}
Let $z$ satisfy \eqref{eq:pig0:extinctS}, then
\begin{align} \label{eq:pig0:lln:hatZ}
\frac{\check{Z}^{(n)}}{n} \convp z \qquad \mbox{as } n \rightarrow \infty.
\end{align}
\end{lemma-pig0-hatZ}
\begin{proof} Using Markov's inequality, \eqref{eq:pig0:lln:hatZ} holds if for $m=1,2$,
\begin{align}  \label{eq:pig0:lln:m}
\E \left[ \left( \frac{\check{Z}^{(n)}}{n} \right)^m \right] \rightarrow z^m \qquad \mbox{as } n \rightarrow \infty.
\end{align}

Since the individuals in the population are exchangeable, we have that
\begin{align*}  
\E \left[  \frac{\check{Z}^{(n)}}{n}  \right] = \P \left( \check{S}_1^{(n)} \neq 0 \right) \rightarrow z \qquad \mbox{as } n \rightarrow \infty.
\end{align*}
Similarly, we have that
\begin{align}  \label{eq:pig0:lln:m2}
\E \left[  \left(\frac{\check{Z}^{(n)}}{n} \right)^2 \right] = \frac{1}{n} \P \left( \check{S}_1^{(n)} \neq 0 \right)  + \frac{n-1}{n} \P \left( \check{S}_1^{(n)} \neq 0,   \check{S}_2^{(n)} \neq 0\right),
\end{align}
where individuals 1 and 2 denote two individuals selected uniformly at random from the population.

Let $E_n$ be the event that two randomly selected individuals 1 and 2 are within a distance $2 \ell_n$ of each other, under the metric $\rho$ defined just after~\eqref{equ:ZnZnch}. Using \eqref{eq:pig0:dist} it is straightforward to show that, for any $\epsilon >0$, the number of individuals within a distance $2 \ell_n$ of individual 1 is at most $n^\epsilon$, for all sufficiently large $n$. Hence, $\P(E_n) \rightarrow 0$ as $n \rightarrow \infty$. Conditional upon $E_n^C$, the growth of the susceptibility sets of individuals $1$ and $2$ over the first $\ell_n$ generations are independent. Therefore
\begin{align}  \label{eq:pig0:lln:En}
 \P \left( \check{S}_1^{(n)} \neq 0,   \check{S}_2^{(n)} \neq 0\right) & =  \P \left( \check{S}_1^{(n)} \neq 0,  \check{S}_2^{(n)} \neq 0, E_n\right)  + \P \left(E_n^C \right) \prod_{j=1}^2  \P \left(\check{S}_j^{(n)} \neq 0|E_n^C \right)  \nonumber \\
 & \rightarrow z^2 \qquad \mbox{as } n \rightarrow \infty.
\end{align}
The lemma follows since \eqref{eq:pig0:lln:En} and \eqref{eq:pig0:lln:m2} imply \eqref{eq:pig0:lln:m} holds for $m=2$.
\end{proof}

In Lemma \ref{lem:pig0:susset} we show, conditional on $Z^{(n)} (= C^{(n)})> \log n$, {\it i.e.}~a major epidemic outbreak has occurred, that  those infected in the epidemic consist mainly of those individuals whose susceptibility set does not go extinct in $\ell_n$ generations.
We prove the lemma under the assumption that $I\equiv 1$, a constant infectious period. In this case the construction of infectious clumps and susceptibility sets are identical, and $\mathbf{V}^{C^{(n)}}_{\ell_n} \eqd \mathbf{V}^{S^{(n)}}_{\ell_n}$, where $\mathbf{V}^{C^{(n)}}_{\ell_n}$ is the number of individuals in the $\ell_n^{{\rm th}}$ generation of the infectious clump. Observe that Lemma~\ref{lem:pig0:susset} yields Theorem~\ref{prop-finsize} with $\pi_G=0$ (recall $R_L=\zeta^S$). In Section~\ref{sec-pig0-genI}, we show how the arguments are adjusted to general $I$.

${\newtheorem{lemma-pig0-susset}[theorem]{Lemma}}$
\begin{lemma-pig0-susset} \label{lem:pig0:susset}
Suppose that $I \equiv 1$ and $\zeta^S > 1$, and let $z$ satisfy \eqref{eq:pig0:extinctS}. Then
\begin{equation}
\label{eq:pig0:susset}
\frac{C^{(n)}}{n}  \;  | \; C^{(n)}> \log n  = \frac{Z^{(n)}}{n} \;  | \;  Z^{(n)} > \log n \convp z, \qquad \mbox{as } n \rightarrow \infty.
\end{equation}
\end{lemma-pig0-susset}
\begin{proof}
Note that $\Cn \eqd \Sn$, since $I \equiv 1$, so~\eqref{eq:pig0:rho2} implies $\lim_{n \to \infty} \P(\Cn > \log n)=z$, where $z>0$ as $\zeta^S > 1$.  Hence, Lemma~\ref{lemma:pig0:hatZ} implies
\begin{equation}
\label{eq:pig0-zhatcond}
\frac{\check{Z}^{(n)}}{n} \;  | \; C^{(n)}> \log n \convp z \qquad \mbox{as } n \rightarrow \infty.
\end{equation}
Thus, to prove \eqref{eq:pig0:susset}, we show that
\begin{align*} 
 \frac{\check{Z}^{(n)}}{n} - \frac{C^{(n)}}{n}  \;  | \;  C^{(n)}> \log n \convp 0 \qquad \mbox{as } n \rightarrow \infty,
 \end{align*}
with the lemma then following immediately from $C^{(n)}= Z^{(n)}$ and~\eqref{eq:pig0-zhatcond}. 

By Markov's inequality, for any $\epsilon >0$,
\begin{align} \label{eq:pig0:Markov}
\P \left( \left|\frac{\check{Z}^{(n)}}{n} - \frac{C^{(n)}}{n} \right| > \epsilon  \;  | \; C^{(n)}> \log n  \right) \leq \frac{1}{\epsilon} \E \left[   \left|\frac{\check{Z}^{(n)}}{n} - \frac{C^{(n)}}{n} \right|   \;  | \;  C^{(n)}> \log n\right].
\end{align}
We can bound the expectation on the right-hand side of \eqref{eq:pig0:Markov} by
\begin{align} \label{eq:pig0:Markov2}
& \E \left[  \left|\frac{\check{Z}^{(n)}}{n} - \frac{C^{(n)}}{n} \right|   \;  | \;  C^{(n)}> \log n\right] \nonumber \\
& \qquad = \E \left[  \frac{1}{n} \left| \sum_{j=1}^n \left( 1_{\{\check{S}_j^{(n)} \neq 0 \}} - 1_{\{ i^\ast \in \mathcal{S}_j^{(n)} \}} \right) \right| \;  | \;   C^{(n)}> \log n\right]  \nonumber \\
& \qquad \leq \frac{1}{n} \sum_{j=1}^n \E \left[  \left| 1_{\{\check{S}_j^{(n)} \neq 0 \}} - 1_{\{ i^\ast \in \mathcal{S}_j^{(n)} \}} \right|   \;  | \;  C^{(n)}> \log n \right] \nonumber \\
& \qquad = \P \left( 1_{\{\check{S}_1^{(n)} \neq 0 \}} \neq 1_{\{ i^\ast \in \mathcal{S}_1^{(n)} \}} \;  | \;   C^{(n)}> \log n \right),
\end{align}
where the final equality in \eqref{eq:pig0:Markov2} follows from the exchangeability of individuals.

Let $\check{E}_n$ be the event that the length of the shortest path between individuals 1 and $i^\ast$ is less than or equal to $\ell_n$. Then $\P (\check{E}_n) \rightarrow 0$ as $n \rightarrow \infty$ and
\begin{align*} 
\P \left(  \check{S}_1^{(n)} = 0, i^\ast \in \mathcal{S}_1^{(n)} \;  | \;  \check{E}_n^C \right) =0.
\end{align*}
Hence, $\P  \left( \check{S}_1^{(n)} = 0, i^\ast \in \mathcal{S}_1^{(n)} | C^{(n)}> \log n \right) \rightarrow 0$ as $n \rightarrow \infty$.

We complete the proof of Lemma \ref{lem:pig0:susset} by showing that
\begin{align} \label{eq:pig0:ss2}
\P  \left( \check{S}_1^{(n)} \neq 0, i^\ast \not\in \mathcal{S}_1^{(n)} | C^{(n)}> \log n \right) \rightarrow 0 \qquad \mbox{as } n \rightarrow \infty.
\end{align}
Following Lemma \ref{lem:pig0:susset1} and the subsequent arguments, for $a=5/(2c^S)$ and $\ell_n = \lceil a \log \log n \rceil$, we have that
\[ \P \left(\left. \mathbf{V}^{C^{(n)}}_{\ell_n}  > [\log n]^2 \right| C^{(n)}> \log n   \right)  \rightarrow 1  \qquad \mbox{as } n \rightarrow \infty. \]

From generation $\ell_n$ we couple the growth of the infectious clump to a lower bound branching process. As $n \to \infty$, the probability that the total number of movers who belong to complexes only through a household (workplace) exceeds $\theta n/2$ tends to 1.   
The number of movers who belong to a given complex only through a household (workplace) is at most $w$. Therefore, conditional on the number of movers exceeding $\theta n/2$, the probability of selecting a given complex, when selecting a mover who belongs to a complex through a household (workplace) at random, is at most $2 w/(\theta n)$. It follows that if an infectious clump is comprised of $K$ complexes then the probability of attempting to add an existing complex from within the clump to the clump is at most $K \times 2 w/(\theta n)$.

Fix $0 < \delta < [\zeta^S -1]/\zeta^S$ and let $\delta^\prime = \delta \theta/(2 w)$. Construct $\mathcal{C}^{(n)}_{i^\ast}$ from generation $\ell_n$ onwards until $\delta^\prime n$ complexes have been added to $\mathcal{C}^{(n)}_{i^\ast}$ as follows. Consider the growth of the infectious clump from each of the $\mathbf{V}^{C^{(n)}}_{\ell_n}  > [\log n]^2$ individuals in generation $\ell_n$ in turn. We choose complexes to add to the infectious clump according to the number of movers in the complex. Whilst fewer than $\delta^\prime n$ complexes have been added to  $\mathcal{C}^{(n)}_{i^\ast}$, the probability we attempt to add an existing complex to the clump is less than $\delta$. We construct a lower bound branching process, $\mathcal{B}^L_\delta$, where the addition of each complex is culled independently with probability $\delta$. Since $[1 -\delta] \zeta^S >1$, the extinction probability, $z_{L,\delta}$, of $\mathcal{B}^L_\delta$ is less than 1. It follows that
\begin{align*} 
\P  \left( C^{(n)}> \delta^\prime n | C^{(n)}> \log n \right) \rightarrow 1 \qquad \mbox{as } n \rightarrow \infty.
\end{align*}

Given that $P ( \check{S}_1^{(n)} >[\log n]^2 |\check{S}_1^{(n)} \neq 0) \rightarrow 1$ as $n \rightarrow \infty$, we have, using the theorem of total probability, that
\begin{align} \label{eq:pig0:ss4}
 \P  \left( i^\ast \not\in \mathcal{S}_1^{(n)} | \check{S}_1^{(n)} \neq 0,C^{(n)}> \log n \right)
& \leq  \P  \left( i^\ast \not\in \mathcal{S}_1^{(n)} | \check{S}_1^{(n)} >[\log n]^2,C^{(n)}> \delta^\prime n \right) \nonumber \\ &\qquad + \P \left(\check{S}_1^{(n)} \leq[\log n]^2| \check{S}_1^{(n)} \neq 0,C^{(n)}> \log n \right) \nonumber \\ & \qquad + \P  \left( C^{(n)}\leq \delta^\prime n | \check{S}_1^{(n)} \neq 0, C^{(n)}> \log n \right).
\end{align}
Using a similar argument to \eqref{eq:pig0:lln:En}, we have that $\P( \check{S}_1^{(n)} \neq 0,C^{(n)}> \log n) \rightarrow z^2$ as $n \rightarrow \infty$. It is then straightforward to show that the latter two terms on the right hand side of \eqref{eq:pig0:ss4} converge to 0 as $n \rightarrow \infty$. Finally
\[\P  \left( i^\ast \not\in \mathcal{S}_1^{(n)} | \check{S}_1^{(n)} >[\log n]^2,C^{(n)}> \delta^\prime n \right) \leq \left[ 1 - \frac{[\log n]^2}{n} \right]^{\delta^\prime n} \rightarrow 0 \qquad \mbox{as } n \rightarrow \infty,
\]
with \eqref{eq:pig0:ss2} following immediately, completing the proof of the lemma.
\end{proof}

\subsection{Non-constant infectious periods} \label{sec-pig0-genI}
In the proof of Lemma \ref{lem:pig0:susset} we have assumed $I \equiv 1$.  We now discuss how the arguments need to be changed for general, random $I$ using the construction of complex infectious clumps given in Section \ref{subsec-outbreakprob}. (Note that the constructions of susceptibility sets do not change and the coupling of $\mathcal{S}^{(n)}$ to $\mathcal{B}^S$ carries over.) In Section \ref{subsec-outbreakprob}, the focus is on the distribution of the severity of an infectious clump whereas here we focus on the probability of extinction of the infectious clump.

Let $\mathbf{Z}^{X,l}$ be defined, as in Section \ref{subsec-outbreakprob}, as a random vector of length $w+h-2$  for the number of movers of each type infected in a complex infectious clump where the initial infective is of type $(X,l)$. 
Let $\bmeta^C = (\eta_{H,1}^C, \ldots, \eta_{H,h-1}^C,\eta_{W,1}^C,\ldots, \eta_{W,w-1}^C)$, where $\eta_{X,l}^C$ is the probability of extinction of the branching process starting from a single complex infectious clump of type $(X,l)$. Then $\bmeta^C$ is the smallest solution in $[0,1]^{w+h-2}$ of
\begin{align*} 
\bmeta^C = \mathbf{g}^C (\bmeta^C) = \left(g^C_{H,1} (\bmeta^C), \ldots, g^C_{H,h-1} (\bmeta^C), g^C_{W,1} (\bmeta^C),\ldots, g^C_{W,w-1} (\bmeta^C) \right)
\end{align*}
where for $\mathbf{x} \in [0,1]^{w+h-2}$,
\begin{align*} 
g^C_{X,l}  (\mathbf{x}) = \E \left[ \prod_{i=1}^{h-1} x_{i}^{Z_{H,i}^{X,l}} \times   \prod_{j=1}^{w-1} x_{j+h-1}^{Z_{W,j}^{X,l}}  \right].
\end{align*}
Then $g^C_{X,l} (\mathbf{x})$, the PGF of $\mathbf{Z}^{(X,l)}$, can be obtained in a similar manner to the PGF of $\mathbf{Z}_X^S$ in Section \ref{subsec-numcomp}. By conditioning on the complex structure given by $M_1, M_2, \ldots, M_d$, defined in Section \ref{subsubsec-clump}, and who the initial infective in the complex infects, $g^C_{X,l} (\mathbf{x})$ can be derived  using Ball and O'Neill \cite{Ball-ONeill-1999}, Theorem 5.1 or Ball \cite{Ball19}, Theorem 4.2.

Returning to the initial infective, let $\mathbf{Z}^R = (Z_{H,1}^R,\ldots,Z_{H,h-1}^R, Z_{W,1}^R, \ldots,Z_{W,w-1}^R)$ be the number of offspring of each type where the initial individual is a remainer with PGF $g^R (\mathbf{x})$. Then we have that
\begin{align} \label{eq:pig0:extinct}
1-\rho = (1-\theta) g^R (\bmeta^C) + \theta \E \left[ g^{H,Q_H} (\bmeta^C) g^{W,Q_W} (\bmeta^C) \right],
\end{align}
with $g^{H,0} (\bmeta^C) = g^{W,0} (\bmeta^C) =1$, since from  $(H,0)$ and $(W,0)$ individuals no infectious clump is created.

By coupling $\mathcal{C}^{(n)}_{i^\ast}$ to the branching process $\mathcal{B}_C$ in analogous fashion to $\mathcal{S}^{(n)}$ to $\mathcal{B}_S$, we have that
\begin{align*}
\P \left( C^{(n)}> \log n\right) \rightarrow \rho \qquad \mbox{as } n \rightarrow \infty,
\end{align*}
where $\rho$ satisfies \eqref{eq:pig0:extinct}.

Let $\tilde{\mathbf{M}}^C$ be the mean offspring matrix of the branching process $\mathcal{B}_C$ (excluding the initial generation) and let $\tilde{\zeta}^C$ be the maximum eigenvalue of $\tilde{\mathbf{M}}^C$. For $0 < \theta <1$, we have that the branching process $\mathcal{B}_C$ is aperiodic and again by Mode \cite{Mode71}, p.19, (8.2), there exists a random variable $W^C$ with $\P (W^C = 0) =1-\rho$ and a continuous distribution on $(0,\infty)$, such that
\begin{align*}
\lim_{k \rightarrow \infty} \frac{1}{[\tilde{\zeta}^C]^k} \mathbf{V}_k^C = W^C \mathbf{v}^C \qquad \mbox{a.s.},
\end{align*}
where $\mathbf{v}^C > \mathbf{0}$ is the normalised left eigenvector of $\tilde{\mathbf{M}}^C$ corresponding to $\tilde{\zeta}^C$. It follows that
\begin{align*} 
\P (\mathbf{V}^{C^{(n)}}_{\ell_n} > [\log n]^2 | C^{(n)} > \log n) \rightarrow 1 \qquad   \mbox{as } n \rightarrow \infty,
\end{align*}
and employing a lower bound branching process as before,
\begin{align}
\label{equ:Cndelta'}
\P (C^{(n)} > \delta^\prime n | C^{(n)} > \log n) \rightarrow 1  \qquad  \mbox{as } n \rightarrow \infty.
\end{align}
As before, the argument is easily modified if $\theta=1$.

We describe how to combine the forward infectious clump with $w+h-2$ types with the backward 2-type susceptibility sets.  From the proof of Lemma \ref{lem:pig0:susset1}, there is, conditional on the susceptibility set surviving to generation $\ell_n$, at least $[\log n]^2$ individuals in generation $\ell_n$ with probability tending to 1 as $n \rightarrow \infty$.  We grow the local infectious  clump $\mathcal{C}^{(n)}_{i^\ast}$ from generation $\ell_n$ in a similar fashion as in the proof of Lemma~\ref{lem:pig0:susset}. However, each time a mover is added to the infectious clump, we determine first who that mover is by sampling uniformly from the set of available movers, including those of the right type ($H$ or $W$) in the $\ell_n^{{\rm th}}$ generation of the susceptibility set $\mathcal{S}_1^{(n)}$.  If the chosen mover belongs to the susceptibility set then individual $1$ is infected by the epidemic and we can stop growing the infectious clump.  Otherwise we determine the finer type of the mover (number of individuals it infects in its new complex) and carry on growing the infectious clump.  In view of~\eqref{equ:Cndelta'},
we can then follow similar arguments to the proof of  Lemma \ref{lem:pig0:susset} to show that an infectious clump containing at least $\delta^\prime n$ individuals will infect a susceptibility set which does not go extinct in the first $\ell_n$ generations with probability tending to 1 as $n \rightarrow \infty$.

Finally, we show that the maximum eigenvalues associated with the forward and backward branching process $\mcB_C$ and $\mcB_S$ are equal, so as is clear on intuitive grounds, $\lim_{n \to \infty} \P(\Cn > \log n)>0$ if and only if $\lim_{n \to \infty} \P(\Sn > \log n)>0$.

${\newtheorem{lemma-zetaS:zetaC}[theorem]{Lemma}}$
\begin{lemma-zetaS:zetaC} \label{lem:zetaS:zetaC}
The maximum eigenvalues of $\tilde{\mathbf{M}}^C$ and $\mathbf{M}^S$ are equal, i.e.
\begin{equation}
\label{eq:zetaS:zetaC}
\tilde{\zeta}^C=\zeta^S.
\end{equation}
\end{lemma-zetaS:zetaC}
\begin{proof}
Consider the spread of a clump within a typical single complex and for $X,Y=H,W$ define $Z^C_{XY}$ analogously to $Z^S_{XY}$.  Recall that in $\mcB_C$, type-$H$ individuals are partitioned into type-$(H,0), (H,1),\dots, (H,h-1)$ individuals and type-$W$ individuals are partitioned into type-$(W,0), (W,1),\dots, (W,w-1)$ individuals.  Note that in the spread of an epidemic within a complex, the infectious period of an individual if it becomes infected is independent of the event that it becomes infected.  Thus, the finer type of each type-$H$ individual infected in the complex epidemic is independently distributed, with probability $p_H(i)$ of being type $(H,i)$, where
\[
p_H(i)=\E\left[\binom{h-1}{i}\left(1-{\rm e}^{-\beta_H' I}\right)^i {\rm e}^{-(h-1-i)\beta_H' I}\right] \quad (i=1,2,\dots h-1),
\]
with $\beta_H'=\beta_H/(h-1)$.  Similarly, each type-$W$ individual is independently type $(W,i)$ with probability $p_W(i)$, where
\[
p_W(i)=\E\left[\binom{w-1}{i}\left(1-{\rm e}^{-\beta_W' I}\right)^i {\rm e}^{-(w-1-i)\beta_W' I}\right] \quad (i=1,2,\dots w-1),
\]
with $\beta_W'=\beta_W/(h-1)$.  (Although the types of individuals are independent, the type of an individual is not independent of the size of the clump unless $I$ is non-random.)

For $i=1,2,\dots, h-1$ and $X=H,W$, let $Z_{(H,i), X}^C$ be the number of type-$X$ offspring of a typical type-$(H,i)$ individual in $\mcB^C$ and define $Z_{(W,i), X}^C$ similarly for $i=1,2,\dots,w-1$.
Exploiting the above-mentioned independence of types, the mean numbers of type-$(H,j_1)$ and type-$(W,j_2)$ individuals in the $k^{\rm th}$ generation of $\mcB^C$ given that the initial individual has type $(H,i)$
are given by
\begin{equation}
\label{equ:meangenkH}
\left(\mu_{Z_{(H,i), H}^C}, \mu_{Z_{(H,i), W}^C}\right) \left(\mathbf{M}^C\right)^{k-1}
\begin{pmatrix}
p_H(j_1) & 0  \\
0& p_W(j_2) \end{pmatrix},
\end{equation}
where
\[
\mathbf{M}^C = \begin{pmatrix} \mu_{Z_{HH}^C} & \mu_{Z_{HW}^C}  \\
 \mu_{Z_{WH}^C}& \mu_{Z_{WW}^C} \end{pmatrix}.
\]
If the initial individual has type $(W,i)$ then these means are given by
\begin{equation}
\label{equ:meangenkW}
\left(\mu_{Z_{(W,i), H}^C}, \mu_{Z_{(W,i), W}^C}\right) \left(\mathbf{M}^C\right)^{k-1}
\begin{pmatrix}
p_H(j_1) & 0  \\
0& p_W(j_2) \end{pmatrix}.
\end{equation}
Note that $\tilde{\zeta}_C$ equals the asymptotic geometric rate of growth of $\E[\mathbf{V}_k^C ]$.  Hence, it follows from~\eqref{equ:meangenkH} and~\eqref{equ:meangenkW} that $\tilde{\zeta}_C=\zeta_C$, where $\zeta_C$ is the maximum eigenvalue of $\mathbf{M}^C$.
To complete the proof we show that the maximum eigenvalues of $\mathbf{M}^C$ and $\mathbf{M}^S$ are equal.

Consider a complex $\mcC$ having structure $\mathbf{A}=\mathbf{a}$, where $a_i$ denotes the number of individuals in group $i$, and at least two movers but without the restriction that group 2 contains at least one mover.  Let $\mcH$ be the set of individuals in $\mcC$ that belong to groups $2,4,\dots 2d$ and $\mcW$ set of individuals in group $2d+1$, so $|\mcH|=|\mcW|=a_{2d+1}$.  Construct on $\mcC$ a random directed graph $\mcG_C$ of potential local contacts in an analogous fashion to $\mathcal{G}^{(n)}$ in Section~\ref{sec-genframework}.  For $i \in \mcH \cup \mcW$, let $\mcC^C_i=\{j \in \mcC: i \leadsto j\}$ and $\mcS^C_i=\{j \in \mcC: j \leadsto i\}$ be respectively the local infectious clump and local susceptibility set of $i$ in $\mcC$.  For $X,Y=\mcH,\mcW$, let $\hat{Z}_{XY}^C(\mathbf{a})=|\mcC^C_i \cap Y|$ and $\hat{Z}_{XY}^S(\mathbf{a})=|\mcS^C_i \cap Y|$, where $i$ is chosen uniformly at random from $X$ and the dependence on the complex structure $\mathbf{a}$ is indicated.  Then, for $X,Y=\mcH,\mcW$,
\begin{align}
\label{equ:clumpsus}
\mu_{\hat{Z}_{XY}^C(\mathbf{a})}=\frac{1}{a_{2d+1}}\sum_{i \in X}\E\left[|\mcC^C_i \cap Y|\right]&=\frac{1}{a_{2d+1}}\E\left[\sum_{i \in X} \sum_{j \in Y} 1_{i \leadsto  j}\right]  \\
&=\frac{1}{a_{2d+1}}
\sum_{j \in Y}\E\left[|\mcS^C_j \cap X|\right]=\mu_{\hat{Z}_{YX}^S(\mathbf{a})}. \nonumber
\end{align}
Observe that~\eqref{equ:clumpsus} continues to hold if for all $i \in \mcH \cup \mcW$, both $\mcC^C_i$ and $\mcS^C_i$ do not include $i$.  With that modification,
\[
\mathbf{M}^C=
\begin{pmatrix}
\E[\mu_{\hat{Z}_{\mcH \mcW}^C(\mathbf{A})}] & \E[\mu_{\hat{Z}_{\mcH \mcH}^C(\mathbf{A})}]  \\
\E[\mu_{\hat{Z}_{\mcW \mcW}^C(\mathbf{A})}]& \E[\mu_{\hat{Z}_{\mcW \mcH}^C(\mathbf{A})}]
\end{pmatrix}
\]
and
\[
\mathbf{M}^S=
\begin{pmatrix}
\E[\mu_{\hat{Z}_{\mcH \mcW}^S(\mathbf{A})}] & \E[\mu_{\hat{Z}_{\mcH \mcH}^S(\mathbf{A})}]  \\
\E[\mu_{\hat{Z}_{\mcW \mcW}^S(\mathbf{A})}]& \E[\mu_{\hat{Z}_{\mcW \mcH}^S(\mathbf{A})}]
\end{pmatrix},
\]
where expectations are with respect to the complex structure $\mathbf{A}$.  Using~\eqref{equ:clumpsus},
\[
\mathbf{M}^S = \begin{pmatrix} \mu_{Z_{HH}^S} & \mu_{Z_{HW}^S}  \\
 \mu_{Z_{WH}^S}& \mu_{Z_{WW}^S} \end{pmatrix}
=\begin{pmatrix} \mu_{Z_{WW}^C} & \mu_{Z_{HW}^C}  \\
 \mu_{Z_{WH}^C}& \mu_{Z_{HH}^C} \end{pmatrix},
\]
so $\mathbf{M}^C$ and $\mathbf{M}^S$ have the same eigenvalues, whence $\zeta^C=\zeta^S$.
\end{proof}

\section{Discussion}\label{sec-disc}

In this paper we studied an epidemic model for a community having two different group structures, households and workplaces, allowed to partly overlap, and asymptotic properties of the epidemic were derived.
For ease of exposition, we have presented our results within the framework of an SIR model in which the infection rates remain constant throughout an infective's infectious period.  Theorems~\ref{prop-pmajor} and~\ref{prop-finsize} can be generalised to the case where infectivity
varies stochastically over an individual's infective period, with $I$ denoting an individual's final cumulative infectivity.

It would of course be interesting to generalise the model to make it more realistic. One obvious generalisation would be to allow households as well as workplaces to have different sizes. Such generalisation is not expected to result in any qualitative new insights. Perhaps more interesting would be to allow individuals to have more heterogeneities: currently the only difference between individuals is their random infectious period. Allowing for varying overall contact rates will have no effect since this could be contained in the infectious period, but allowing how the contact rate divides between household, workplace and community to be different, and also allowing susceptibility to vary, would be an interesting as well as from an applied point of view relevant, extension of the model. Another interesting extension would be to consider non-pharmaceutical interventions (NPIs), typically acting differently on the three types of contact, to see how this alters an ongoing epidemic.

A natural theoretical extension would be to supplement the law of large numbers for the size of major outbreak in Theorem~\ref{prop-finsize} with a central limit theorem.  Figure~\ref{fig:finalsize1} in Section~\ref{sec-illustrations} suggests that $\Zn|\Zn>\log n$ is asymptotically normal when $\pi_G>0$ and the embedding construction in Section~\ref{subsec-finsizeproof} offers a possible method of proof (cf.~Ball et al.~\cite{BMST} for the households model).  However, except in the special case $\theta=0$, proving a functional central limit theorem analogue of Lemma~\ref{lem:glivenko} is not straightforward and a useful expression for the asymptotic variance is difficult to obtain, owing to overlap of complexes.

The numerical illustrations suggest that the final size $z$ increases with  $\theta$ (the amount of non-overlap between group structures) and workplace size $w$.  Other numerical calculations (not shown) suggest that $z$ also increases with $\pi_G$ (the fraction of contacts that are global) and household size $h$.  An obvious open problem is to prove these results in generality, or to find counterexamples.

Another direction for future work is investigation of vaccination strategies under our model.  Patwardhan et al.~\cite{P23} use Monte Carlo simulations to show that when a fraction of the population is vaccinated with a perfect vaccine (i.e.~one which necessarily renders the recipient fully immune) by selecting individuals for vaccination uniformly at random, there is a greater reduction in both the peak and final size of an outbreak when there is full overlap ($\theta=0$) than when there is no overlap ($\theta=1$), provided $\beta$ is calibrated suitably for the two models.  The methodology of this paper can be extended to allow the effect of such vaccination on final size to be investigated more systematically.  Other vaccination-related questions could also be explored.  For example, is it better to target vaccination at movers or remainers, and what is the optimal vaccine allocation strategy for a given vaccination coverage?

Although the reproduction numbers $R_*$ and $R_L$ serve as threshold parameters for the model with $\pi_G>0$ and $\pi_G=0$, respectively, they do have some limitations; see Pellis et al.~\cite{Pellis09} for a discussion of $R_*$ for the model with no overlap ($\theta=1$).  For example, the fact that $R_*$ can be infinite means that it can be completely uninformative about the effort that is required to bring an epidemic below threshold.   Even though $R_L$ is necessarily finite, it can be misleading when comparing different models (for example, models with different values of $\theta$), since it is based on the proliferation of infected complexes, rather than infected individuals, and different models may have different complex size distributions.  A similar comment applies to $R_*$.  It is possible to generalise the individual-based reproduction numbers $R_0$ and $R_I$, given for the case $\theta=1$ in Ball et al.~\cite{BPT16}, Sections 4.2 and 4.5, respectively, to general $\theta$, though numerical calculation of the former is generally prohibitive.  These will be considered in a subsequent paper.

\section*{Acknowledgements}
Tom Britton is grateful to the Swedish Research Council (grant 2020-0474) for financial support.

\end{document}